\documentclass[11pt,a4paper,reqno]{amsart}
\usepackage{amssymb,amsmath,amsthm}
\usepackage[utf8]{inputenc}
\usepackage[T1]{fontenc}
\usepackage{enumerate}
\usepackage[all]{xy}
\usepackage{mathrsfs} 
\usepackage{fullpage}
\usepackage{comment}
\usepackage{xcolor}

\def\Z{\mathbb{Z}}
\def\P{\mathbb{P}}

\def\cA{\mathcal{A}}
\def\sL{\mathscr{L}}
\def\cM{\mathcal{M}}
\def\cA{\mathcal{A}}

\def\sC{\mathscr{E}}

\def\cO{\mathcal{O}}

\def\F{\mathbb{F}}
\def\N{\mathbb{N}}
\def\Q{\mathbb{Q}}

\def\div{\mathrm{div}}
\def\ord{\mathrm{ord}}

\newtheorem{thm}{Theorem}[section]
\newtheorem{lem}[thm]{Lemma}

\newtheorem{cor}[thm]{Corollary}

\newtheorem{claim}{Claim}
\theoremstyle{definition}
\newtheorem{defn}[thm]{Definition}

\theoremstyle{remark}
\newtheorem{rmk}[thm]{Remark}
\title{Counting integral points of fixed degree and given height over function fields}

\author{Angelot Behajaina}
\email{behajainaangelot@gmail.com}
\address{Department of Mathematics, Institute of Science Tokyo, 2-12-1 Ookayama, Meguro-ku, Tokyo 152-8551, Japan}

\author{Susheel Shankar}
\email{susheel.shankar@tugraz.at}
\address{TU Graz, Institute of Analysis and Number Theory, Steyrergasse 30/II, 8010 Graz, Austria}

\begin{document}

\maketitle
\begin{abstract}
Let $k$ be a finite extension of the rational function field in one variable over a finite field, and let $\cO_k$ be its ring of integers. Let $\overline{k}$ be an algebraic closure of $k$. For integers $e \geq 2$ and $n \geq 1$, we study the number of points in the affine space $\overline{k}^n$ whose coordinates are integral over $\cO_k$, have absolute logarithmic Weil height equal to $m \geq 0$, and generate an extension of degree $e$ over $k$. More precisely, we derive an asymptotic formula for this quantity as $m \rightarrow \infty$ whenever $n > e$, except in the cases $(n,e)\in \{(3,2),(4,3),(5,4)\}$. When $k$ is the rational function field, however, we use Weierstrass semigroups to obtain asymptotic formulas for the cases $(n,e)\in \{(1,2),(2,2),(3,2)\}$, thereby completing the picture in the setting $e=2$. As an application, we also derive asymptotic estimates for the number of certain decomposable forms with coefficients in $\cO_k$.
\end{abstract}

\section{Introduction}
In this article, we count algebraic points with integral coordinates, prescribed degree, and fixed Weil height over a function field.

The corresponding problem over number fields has been extensively studied. Let $k$ be a number field of degree $d$ over $\Q$, let $\overline{k}$ be a fixed algebraic closure of $k$, and let $n \geq 1$. By Northcott's theorem \cite{Nor49}, the set of points in $\mathbb{P}^{n}(\overline{k})$ of bounded degree and bounded multiplicative Weil height $H$ is finite. It is therefore natural to seek asymptotic estimates as the height bound tends to infinity. Schanuel \cite{Sch79} proved that
$$
N_k(n,1,X):=\#\{ P \in \mathbb{P}^{n}(k) \mid H(P) \leq X \}=S_k(n)X^{d(n+1)}+O(X^{d(n+1)-1}\log(X)),\quad\textrm{as }X \rightarrow \infty,
$$ where $S_k(n)$ is a positive constant depending only on $k$ and $n$, commonly called \emph{Schanuel's constant}. More generally, for $e \geq 1$, Schmidt \cite{Sch93} considered
$$
N_k(n,e,X):=\#\{ P=(\alpha_0:\dots:\alpha_n) \in \mathbb{P}^{n}(\overline{k}) \mid [k(P):k]=e, H(P) \leq X \},
$$ where $k(P)$ is the field generated over $k$ by the ratios $\alpha_i/\alpha_j$. He established bounds of the form
$$
 c(k,e,n)X^{d e(\max\{e,n\}+1)}\leq N_k(n,e,X) \leq C(k,e,n) X^{d e(e+n)},
$$ for suitable positive constants $c(k,e,n)$ and $C(k,e,n)$. Over $\Q$, Schmidt \cite{Sch95} later obtained an asymptotic formula for the case $e=2$, while Gao \cite{Gao95} established one for $n>e>2$. Over an arbitrary number field, Masser and Vaaler \cite{MV07} proved the full asymptotic when $n=1$, and Widmer \cite{Wid09} obtained an asymptotic formula under the condition $n > 5 e/2 +4 +2/(d e)$. 

Analogous questions over function fields have also received considerable attention. Thunder and Widmer \cite{TW13} proved asymptotic estimates for $n> 2e+2$, while Kettlestrings and Thunder \cite{Thu15} treated the quadratic case $e=2$ with $n \geq 2$. More recently, Thunder \cite{Thu25} established asymptotic formulas throughout the range $n \geq e+1$, together with sharp order estimates for $1 \leq n \leq e$. 

A related problem concerns algebraic points with integral coordinates in affine space. Given a number field $k$ with ring of integers $\cO_k$, one considers 
$$
N_{\cO_k}(n,e,X)=\#\{ P \in \overline{\cO_k}^n \mid [k(P):k]=e, H(P) \leq X \},
$$ where $\overline{\cO_k}$ denotes the integral closure of $\cO_k$ in $\overline{k}$. Results in this direction go back to Lang \cite[p. 81]{Lan83}, who treated $(n,e)=(1,1)$. Chambert-Loir and Tschinkel \cite[Theorem 3.5.6]{CT12} generalized this to $n\ge 1, e=1$.  Chern and Vaaler \cite{CV01} handled the case $k=\mathbb Q, n=1, e\ge 1$, which was generalized by Barroero\cite{Bar14} to arbitrary $k$ and $n=1$. Widmer \cite{Wid16} obtained asymptotic formulas in the range $n \geq e+8$ for an arbitrary ground field $k$.

By contrast, the corresponding problem over function fields does not appear to have been investigated previously. The purpose of the present paper is to initiate such a study. 

Let $k/\F_q(t)$ be a finite extension of degree $d=[k:\F_q(t)]$, fix an algebraic closure $\overline{k}$, and assume that $k \cap \overline{\F_q}=\F_q$. For every finite extension $K/k$, set $q_K:=\#(K \cap \overline{\F_q})$, and let $\cO_K$ be the integral closure of $\F_q[t]$ in $K$. We write $\overline{\cO_k}$ for the integral closure of $\cO_k$ in $\overline{k}$. 

Throughout the paper, we use the absolute logarithmic Weil height $h$; see Section \ref{ss:height}. For $P=(\alpha_1,\dots,\alpha_n) \in \overline{k}^n$, write $k(P):=k(\alpha_1,\dots,\alpha_n)$ for the field generated over $k$ by its coordinates. For integers $n \geq 1, e\geq 1$, and $m \geq 0$, define 
$$
N_{\cO_k}(n,e,m):=\#\left\lbrace P \in \overline{\cO_k}^n \mid [k(P):k]=e, h(P)=\frac{m}{ed}, q_{k(P)}=q \right\rbrace.
$$ 

We extend this notation by setting $N_{\cO_k}(n,d,m)=0$ whenever $d \text{ or } m$ is not an integer.

Our principal aim is to determine the asymptotic behaviour of $N_{\cO_k}(n,e,m)$ as $m \rightarrow \infty$.

Let $K/k$ be a finite extension of degree $e$ with $q_K=q$. Write $g_K$ for its genus, and let $\cM^\infty_K$ denote the set of places of $K$ lying above the place $\infty$ of $\F_q(t)$. Let $\mathbb N$ denote the set of non-negative integers.
Set $$
U_K=\sum_{v \in \cM^\infty_K} \N \cdot \deg(v).
$$ We shall see in Section \ref{ss:height} that the height $h$ of a point counted over $K$ takes values in $U_K/(ed)$. 

Define 
$$
s_K:=\#\cM^\infty_K,\quad \rho_K:=\gcd\{\deg(v) \mid v \in \cM^\infty_K \},
$$ and 
$$
\Theta_{n,K}:=\prod_{v \in \cM^\infty_K} \left(1 - \frac{1}{q^{n \deg(v)}} \right).
$$ We further define
$$
\kappa_n(K):=\frac{\rho_K}{(s_K-1)! \prod_{v \in \cM^\infty_K} \deg(v)} \cdot \Theta_{n,K}\, ,
$$
and
$$
S_{\cO_K}(n)=\kappa_n(K)q^{-ng_K+n},
$$
which is the integral analogue of Schanuel's constant. Finally, set
$$
S_{\cO_k}(n,e)=\sum_{K}S_{\cO_K}(n),
$$
where the sum runs over all degree-$e$ extensions $K/k$ with $q_K=q$ for which the places in $\cM_k^\infty$ splits completely in $K$. Equivalently, $s_K=es_k$, the largest possible value in a degree-$e$ extension.

\begin{thm}\label{thm:main1}
 Assume that $e \geq 2$ and $n \geq e+1$.  The series $S_{\cO_k}(n,e)$ converges and as $m \rightarrow \infty$ through $m \in U_k$, we have for  
$
(n,e) \notin \{(3,2),(4,3),(5,4)\},
$ 
\begin{align*}
N_{\cO_k}(n,e,m)&=S_{\cO_k}(n,e) m^{s_ke-1}q^{mn}+O_{k,n,e} \left(m^{s_ke-1}q^{\frac{(n+3)m}{2}}+m^{s_ke-2}\, q^{mn}\right).
\end{align*}

When $(n,e)=(3,2),$ we have
$$
N_{\cO_k}(3,2,m)=
O_{k,e}\left(q^{\frac{7m}{2}}\right).
$$

Moreover, we have for all $\varepsilon>0,$
\begin{align*}
N_{\cO_k}(e+1,e,m)&=
 \begin{cases}
O_{k,e, \varepsilon}\left(q^{4m+m\varepsilon}\right)&\textrm{if }e= 3,\\
O_{k,e,\varepsilon}\left(q^{5m+\frac{m}{3}+m\varepsilon}\right)&\textrm{if }e=4.
\end{cases}
\end{align*}
The implied constants depend only on $k$, $n$ and $e$, except for the cases $(n,e)=(4,3)\text{ and }  (5,4),$ where it additionally depends on $\varepsilon$.
\end{thm}
\noindent

In the special case when $k=\F_q(t)$, we use the theory of  Weierstrass semigroups to provide the asymptotics in the cases $(n,e)\in \{(3,2),(2,2),(1,2)\}$.
\begin{thm}\label{thm:main2}
     If $k=\F_q(t)$ and $q$ is odd, as $m \rightarrow \infty$,
$$
N_{\cO_k}(1,2,m)={2\left(q^2-1\right)}q^{2m}+O(mq^{m}),
$$
$$
N_{\cO_k}(2,2,m)={\frac{3}{4}q^3(q-1)\left(1-\frac{1}{q^2} \right)^2} m^2 q^{2m}+O(m q^{2m}),
$$ and
$$
N_{\cO_k}(3,2,m)= S_{\cO_{k}}(3,2)mq^{3m}+O(q^{3m}).
$$
\end{thm}

The proof of Theorem \ref{thm:main1} and Theorem \ref{thm:main2} rely on counting integral points of fixed height in a given function field. Let $K/\F_q(t)$ be a finite extension of degree $\ell$ with $q_K=q$, and set
$$
D_K=\sum_{v \in \cM^\infty_K}\deg(v).
$$

\begin{thm}\label{thm:poinbounheiK}
For $m \in U_K$ with $m>0$,
$$
N_{\cO_K}(n,m)
:= \#\left\lbrace \alpha \in \cO_K^n \mid h(\alpha)=\frac{m}{\ell} \right\rbrace
=S_{\cO_K}(n) m^{s_K-1}\, q^{mn}
+ O_{s_K}(E(m)),
$$
where, for $s_K \geq 2$,
$$
E(m)=
\begin{cases}
m^{s_K-2}\, q^{mn-ng_K+n}
& \textrm{if } m \geq 2 g_K-1+D_K,\\[0.5em]
\displaystyle
m^{s_K-1} q^{\frac{m}{2}n+n}+m^{s_K-2}q^{mn-ng_K+n}
& \textrm{if }  2g_K-1 \leq m \leq 2g_K-2+ D_K,\\[0.5em]
\displaystyle
m^{s_K-1}q^{\frac{m}{2}n+n}+ m^{s_K-2}q^{n g_K-n} 
& \textrm{if } m < 2g_K-1,
\end{cases}
$$
and, for $s_K=1$,
$$
E(m)=
\begin{cases}
0 & \textrm{if } m \geq 2 g_K-1+D_K,\\
q^{\frac{m}{2}n+n} & \textrm{otherwise}.
\end{cases}
$$
The implied constant depends only on $s_K$.
\end{thm}

We now assume that $K/k$ has degree $e$ and satisfies $q_K=q$. The proof of Theorem \ref{thm:main1} uses the following consequence of Theorem \ref{thm:poinbounheiK}. Following \cite{TW13}, define
$$
\delta_{n+1}(K/k)=\min\{h(P)\mid P \in \mathbb{P}^{n}(K)\,\,\textrm{and}\,\, k(P)=K\}.
$$

For $m \in \N$, define
$$
N_{\cO_k}(n,K,m)
:= \#\left\lbrace \alpha \in \cO_K^n \mid h(\alpha)=\frac{m}{de}, k(\alpha)=K \right\rbrace,
$$
and set $N_{\cO_k}(n,K,m)=0$ whenever $m \notin \N$.

\begin{cor}\label{cor:ofdegr}
 Let $m \in U_K$ with $m>0$, and suppose that $\delta_{n+1}(K/k) \leq m/(de)$. Then 
$$
N_{\cO_k}(n,K,m)= S_{\cO_K}(n) m^{s_K-1}\, q^{mn}
+ O_{k,n,e}(\tilde E(m)),
$$
where, for $s_K \geq 2$, 
$$
\tilde E(m)=
\begin{cases}
m^{s_K-2}\, q^{mn-ng_K}+q^{\frac{m(n+1)}{2}}
& \textrm{if } m \geq 2 g_K-1,\\[0.5em]
\displaystyle
m^{s_K-2}q^{n g_K-n}+q^{\frac{m(n+1)}{2}}
& \textrm{otherwise},
\end{cases}
$$
and, for $s_K=1$,
$$
\tilde E(m)=q^{\frac{m(n+1)}{2}}.
$$
\end{cor}

As an application, we count irreducible decomposable forms with integral coefficients. This is an affine analogue of the projective counting problem studied in \cite{TW13}.

Let $e \geq 1$ and let
$$
F= \sum_{\substack{I=(i_0,\dots,i_n) \in \N^{n+1}\\ i_0+\dots+i_n=e}} c_I{\bf x}^I \in \cO_k[x_0,\dots,x_n],\quad {\bf x}^{I}=x_0^{i_0}\cdots x_n^{i_n}, 
$$
be a homogeneous form of degree $e$, irreducible over $k$, and monic in the sense that 
$$
c_{(e,0,\dots,0)}=1.
$$
The height of $F$, denoted by $h(F)$, is defined to be the absolute logarithmic Weil height of the coefficient vector of $F$, viewed as a point of $\P^{\binom{n+e}{e}-1}(k)$. We say that $F$ is \emph{decomposable} of \emph{type $(n,e)$} if, over $\overline{k}$, it admits a factorization 
\begin{equation}\label{eq:decomlin}
F=\prod\limits_{j=1}^e L_j,\quad L_j=x_0+\alpha_{j,1}x_1+ \dots+ \alpha_{j,n}x_n,
,\end{equation} 
and if, for one (equivalently, every) factor,
$$
K_j:=k(\alpha_{j,1},\dots,\alpha_{j,n})\quad\textrm{satisfies } [K_j:k]=e\textrm{ and }q_{K_j}=q. 
$$
Let 
$N_{\cO_k}(n,e,m,\textrm{hom})$ denote the number of such forms satisfying
$h(F)=m/ed$. We set this number to be zero if $e$ or $m$ is not an integer.

Let $p={\rm char}(k)$, and write $e=p^r e_0$ with $p \nmid e_0$. For $i=0,1, \dots,r$, define
 $$
 T_i=N_{\cO_k}\left(n,\frac{e}{p^i},\frac{m}{p^i}\right)+(p-1)\sum\limits_{j=i+1}^{r}p^{j-i-1}N_{\cO_k}\left(n,\frac{e}{p^j},\frac{m}{p^i}\right).
 $$

\begin{thm}\label{hom_forms_gen}
 For every $e \geq 1$ and $m \geq 0$, we have 
   \begin{equation}
       e N_{\cO_k}(n,e,m,{\rm hom})=T_0-(p-1)\sum\limits_{i=1}^{r}T_i.    
   \end{equation}
\end{thm}

\begin{cor}\label{cor: hom_forms_gen} Let $e \geq 2$ and $n \geq e+1$, with
$$
(n,e)\notin\{(3,2),(4,3),(5,4)\}.
$$
As $m \rightarrow \infty$ through $ m \in U_k$, 
    \begin{align*}
        e N_{\cO_k}(n,e ,m,{\rm hom})= S_{\cO_k}(n,e) m^{s_k e-1}q^{mn}+O_{k,n,e} \left(m^{s_k e-2}\, q^{mn}\right).
    \end{align*}
Furthermore, if $k=\F_q(t)$ and $q$ is odd, as $m \rightarrow \infty$,
    \begin{align*}
         &N_{\cO_k}(1,2,m,{\rm hom})=\left(q^2-1\right)q^{2m}+O(mq^{m}),\\
&N_{\cO_k}(2,2,m,{\rm hom})={\frac{3}{8}q^3(q-1)\left(1-\frac{1}{q^2} \right)^2} m^2 q^{2m}+O(m q^{2m}),\\
&N_{\cO_k}(3,2,m,{\rm hom})= \frac 12 S_{\cO_{k}}(3,2)mq^{3m}+O(q^{3m}),
    \end{align*}
    
\end{cor}

\subsection{Organization of the paper} Section \ref{ss:height} reviews the necessary background on function fields, heights, Riemann-Roch spaces, and quadratic extensions. Section \ref{sec:poibouhei} proves Theorem \ref{thm:poinbounheiK}. In Section \ref{ss:prfmain1}, we deduce Corollary \ref{cor:ofdegr} and prove Theorem \ref{thm:main1}. Section \ref{sec:smalcas} proves Theorem \ref{thm:main2}. Finally, Section \ref{sec:last} proves Theorem \ref{hom_forms_gen} and Corollary \ref{cor: hom_forms_gen}.

\section{Preliminaries on Function Fields}\label{ssec:basicff}\label{ss:height}
In this section, we recall the facts about function fields, heights, Riemann-Roch spaces, and quadratic extensions that will be used throughout the paper.

Let $K/\F_q(t)$ be a finite extension of degree $\ell$ with $q_K=q$. Let $\mathcal{M}_K$ denote the set of places of $K$. For each $v \in \mathcal{M}_K$, let $K_v$ be the completion of $K$ at $v$, and let 
$$
\ord_v:K_v \rightarrow \Z \cup \{\infty\}
$$ be the normalized valuation, so that $\ord_v(K_v)=\Z$. For $y=(y_0:\dots:y_n)\in \P^n(K)$, define 
$$
\ord_v(y)=\min_{0 \leq i \leq n}\ord_v(y_i).
$$
Similarly, for $\alpha=(\alpha_1,\dots,\alpha_n)\in K^n$, define
$$
\ord_v(\alpha)=\min_{1 \leq i \leq n}\ord_v(\alpha_i).
$$
The \emph{degree} of a place $v$ is the degree of its residue field over $\F_q$. 

\begin{defn}\label{defn:relweihei}
For $y=(y_0:\dots:y_n)\in \P^n(K)$, the \emph{relative (logarithmic) Weil} height is $$h_K(y)=-\sum\limits_{v \in \cM_K}\ord_\nu(y)\cdot \deg(v),$$
and the \emph{absolute (logarithmic) Weil} height is 
$$h(y)=\frac{h_K(y)}{\ell}.$$ 
For $\alpha \in K^n$, we set
$$
h_K(\alpha):=h_K(1:\alpha_1\dots:\alpha_n)\quad\textrm{and}\quad h(\alpha):=h(1:\alpha_1:\dots:\alpha_n).
$$
\end{defn}
\begin{rmk}
Provided that $q_K=q$, the height $h(y)$ (resp., $h(\alpha)$) is independent of the choice of the field $K$ containing $y$ (resp., $\alpha$).
\end{rmk}

Let 
$$
F=\sum\limits_{I}c_I{\bf{x}}^I \in K[x_0,\dots,x_n]
$$ be a nonzero polynomial of total degree $\delta$. We define $h(F)$ to be the absolute logarithmic Weil height of its coefficient vector, viewed as a point of
$$
\P^{^{\binom{n+1+\delta}{\delta}-1}}(K).
$$

The following lemma relates the height of a polynomial to the heights of its factors.
\begin{lem}\label{lem:heighaddpoly}
Let $P,P_1,P_2 \in K[x_0,\dots,x_r]$ be nonzero polynomials. If $P=P_1P_2$, then
$$
h(P)=h(P_1)+h(P_2).
$$
\end{lem}

\begin{proof}
Let $A,A_1,A_2$ be the coefficient vectors of $P,P_1,P_2$, respectively. By Gauss's Lemma \cite[Lemma 1.6.3]{BG06}, for every $v \in \cM_K$, 
$$
\ord_v(A)=\ord_v(A_1)+\ord_v(A_2).
$$
Summing over all places gives
$$
h(P)=h(A)=h(A_1)+h(A_2)=h(P_1)+h(P_2).
$$ 
\end{proof}

\begin{rmk}
Although \cite[Lemma 1.6.3]{BG06} is stated for number fields, its proof applies to every non-archimedean place; see also \cite[\S1.6.4]{BG06}. 
\end{rmk}

Let ${\rm Div}(K)$ denote the divisor group of $K$. Thus, a divisor is a finite formal sum 
$$A=\sum_{v \in \mathcal{M}_K}A_v \cdot v, \quad A_v \in \Z.$$ The \emph{support } of $A$ is
$$
{\rm supp}(A)=\{ v \in \cM_K \mid A_v \neq 0\}.$$ 
We write 
$$A_0=\sum_{v \in \cM_K, A_v>0} A_v \cdot v \quad \textrm{and} \quad A_\infty=\sum_{v \in \cM_K, A_v<0} -A_v \cdot v
$$ for the \emph{zero}  and \emph{pole} parts of $A$, respectively, so that $A=A_0-A_\infty$.

It is well known that every nonzero $\alpha \in K^n$ determines a divisor  
$$
(\alpha)=\div(\alpha)=\sum_{v \in \mathcal{M}_K} \ord_v(\alpha) \cdot v.
$$

We have the following lemma, that follows immediately from the definitions:
\begin{lem}
For every $\alpha \in K^n$,
$$
h_K(\alpha)=\deg((\alpha)_\infty)=-\deg(\div(1,\alpha))=\sum_{v \in \mathcal{M}_K}\max\{0,-\ord_v(\alpha)\} \deg(v).
$$
\end{lem}

\begin{rmk}\label{rmk:heihpropre} Let $\cM_K^\infty$ denote the set of places of $K$ lying above $\infty$. Then
$$
\cO_K=\{\alpha \in K \mid \ord_v(\alpha) \geq 0\,\,\textrm{for all }v \in \mathcal{M}_K\setminus \cM^\infty_K\}.
$$ Consequently 
$$h_K(\cO_K) \subset U_K \quad \textrm{and}\quad U_K=\sum_{v \in \cM_K^\infty} \N \cdot \deg(v).
$$
\end{rmk}

\subsection{Riemann--Roch spaces}
For $A \in {\rm Div}(K)$, the \emph{Riemann--Roch space associated with $A$} is the $\F_q$-vector space
$$
\mathscr{L}_n(A):=\{ \alpha \in K^n\mid (\alpha)+A \geq 0 \}=\{\alpha \in K^n \mid \ord_v(\alpha) \geq -A_v\,\,\textrm{for all } v \in \cM_K\}.
$$ Set $\ell_n(A):=\dim_{\F_q}(\mathscr{L}_n(A))$. Since
$$
\mathscr{L}_n(A)=(\mathscr{L}_1(A))^n,
$$
we have 
\begin{equation}\label{rieroch_n}
\ell_n(A)=n\ell_1(A).    
\end{equation}

We shall use the following standard facts; see, for example, \cite{Sti09}.
\vskip 1mm
\noindent
$\bullet$ If $\deg(A)<0$, then $\ell_1(A)=0$.
\vskip 1mm
\noindent
$\bullet$ ({\bf Clifford's theorem}) If $0 \leq \deg(A) < 2 g_K-1$, then
$$
\ell_1(A) \leq 1 + \frac{1}{2} \deg(A).
$$
\vskip 1mm
\noindent
$\bullet$ ({\bf Riemann--Roch theorem}) If $\deg(A) \geq 2g_K-1$, then
$$
\ell_1(A)=\deg(A)-g_K+1. 
$$
We shall also use the following estimate; see \cite[Lemma 1.4.8]{Sti09}.
\begin{lem}\label{lem:bounddima}
Let $A, B \in {\rm Div}(K)$ with $A \leq B$. Then 
$$
\dim_{\F_q}(\mathscr{L}_1(B))- \dim_{\F_q}(\mathscr{L}_1(A)) \leq \deg(B)-\deg(A).
$$
\end{lem}

\subsection{Quadratic function field} In this subsection, we collect several facts about quadratic extensions of 
$$
k_0:=\F_q(t).
$$  
Throughout the subsection, we assume that $q$ is odd. 

\subsubsection{Basic properties} 

\begin{lem}\label{lem:sevres} Let $K=k_0(\sqrt{f})$, where $f \in \F_q[t]$ is a nonconstant squarefree polynomial. Then the following statements hold.
\begin{enumerate}
\item The full constant field of $K$ is $\F_q$, and
$$
g_K=
\begin{cases}
(\deg(f)-1)/2 &\textrm{if } \deg(f) \textrm{ is odd},\\ 
(\deg(f)-2)/2 &\textrm{if } \deg(f) \textrm{ is even}.
\end{cases}
$$
\item Let $\lambda$ be the leading coefficient of $f$.
\begin{itemize}
    \item If $\deg(f)$ is odd, then $\infty$ ramifies in $K/k_0$.
    \item If $\deg(f)$ is even and $\lambda$ is a square in $\F_q$, then $\infty$ splits in $K/k_0$.
    \item If $\deg(f)$ is even and $\lambda$ is not a square in $\F_q$, then $\infty$ is inert in $K/k_0$.
\end{itemize}
\end{enumerate}
\end{lem}

\begin{proof}
Part (1) follows from \cite[Example 3.7.6.]{Sti09}. Part (2) is \cite[Proposition 14.6]{Ros02}.
\end{proof}
The preceding lemma yields the following count.
\begin{lem}\label{lem:gen_count}
For every integer $i \geq 0$,
$$
\sum\limits_{\substack{[K:k_0]=2\\q_K=q\\g_K=i\\\infty \textrm{ ramifies in }K/k_0}}1=
\begin{cases}
2q &\textrm{ if } i=0,\\
2q^{2i}(q-1) &\textrm{ if } i>0.
\end{cases}
$$
$$
\sum\limits_{\substack{[K:k_0]=2\\q_K=q\\g_K=i\\\infty \textrm{ inert in }K/k_0}}1= q^{2i+1}(q-1),
$$
and
$$
\sum\limits_{\substack{[K:k_0]=2\\q_K=q\\g_K=i\\\infty \textrm{ splits in }K/k_0}}1= q^{2i+1}(q-1).
$$
Consequently, $$\sum\limits_{\substack{[K:k]=2\\q_K=q\\g_K=i}} 1=
\begin{cases}
2q^2 &\textrm{ if } i=0, \\
2(q^2-1)q^{2i} &\textrm{ if } i>0.
\end{cases}
$$
\end{lem}
\begin{proof} 
From \cite[Proposition 2.3]{Ros02} we have, 
\begin{align*}
\#\{ f\in \F_q[t] \mid \deg(f)=n,\,f\textrm{ is monic and  square-free}\}
=\begin{cases}
q&\textrm{if }n=1,\\
(q^n-q^{n-1})&\textrm{if }n\geq 2.
\end{cases}
\end{align*}

And therefore, it follows that
\begin{align*}
\#\{ f\in \F_q[t] \mid \deg(f)=n,\,f\textrm{ is square-free}\}
=\begin{cases}
q(q-1)&\textrm{if }n=1,\\
(q-1)(q^n-q^{n-1})&\textrm{if }n\geq 2.
\end{cases}
\end{align*}

The result now follows from Lemma \ref{lem:sevres} and the fact that, for nonconstant squarefree $f,g \in \F_q[t]$,  $$
k_0(\sqrt{f})=k_0(\sqrt{g}) \text{ if and only if } f/g \in (\F_q^\times)^2.
$$
\end{proof}

\subsubsection{Weierstrass semigroups}\label{ssec:weismg} For quadratic extensions of $k_0$, heights can be described explicitly in terms of the possible pole orders at the places above $\infty$. Given $K/k_0$ with $q_K=q$, and places $P_1,P_2\in \cM_K^{\infty}$ , we define the \emph{Weierstrass semigroups} 
\begin{align*}
W(P_1)&:=\{ r \in \N \mid (\alpha)_\infty=r P_1\textrm{ for some }\alpha \in K\}
\end{align*} and
\begin{align*}
W(P_1,P_2)&:=\{(r_1,r_2) \in\N^2 \mid (\alpha)_\infty=r_1P_1+r_2P_2\,\textrm{ for some } \alpha \in K\}.
\end{align*}

\begin{lem}\label{lem:WeieGap} Let $K/k_0$ be a quadratic extension with $q_K=q$.
\begin{enumerate}
\item Let $P$ be a place of $K$ above $\infty$. Then
\begin{align}
W(P)=\begin{cases}
\N &\textrm{if }\infty\textrm{ is inert},\\
\N \setminus\{1,3,\dots,2g_K-1\}&\textrm{if }\infty\textrm{ is ramified},\\
\N \setminus\{1,2,\dots,g_K\}&\textrm{if }\infty\textrm{ splits}.
\end{cases}
\end{align}
\item Suppose that $\infty$ splits into two places $P_1$ and $P_2$ of $K$. Then
\begin{align*}
W(P_1,P_2)=\N^2 \setminus \{(i,j) \mid 0 \leq i\neq j \leq g_K\}.
\end{align*}
\end{enumerate}
\end{lem}

\begin{proof} When $k_0$ is the rational function field over an algebraically closed field, part (1) follows from \cite[Exercise 6.12]{Sti09}, while part (2) follows from \cite[Example 3.2]{Hom96}. For completeness, we give a direct proof of part (1).

Suppose first that $\infty$ is inert in $K$. Then 
$$
\ord_P(t^i)=i\cdot \ord_{\infty}(t)=-i
$$ for every $i \geq 0$, and hence $W(P)=\N$.

Suppose next that $\infty$ is ramified in $K$. Then 
$$
\ord_P(t^i)=2i\cdot \ord_{\infty}(t)=-2i
$$ for all $0 \leq i \leq g_{K}$. Together with \cite[Theorem 1.6.8]{Sti09}, this gives 
$$
W(P)=\N \setminus\{1,3,\dots,2g_K-1\}.
$$

Finally, suppose that $\infty$ splits into two places $P_1$ and $P_2$ of $K$. For $r<g_K$, then Clifford's theorem gives
$$
\ell_1(r(P_1+P_2))\le 1+r.
$$
On the other hand, Riemann-Roch Theorem gives
$$\ell_1(g_K(P_1+P_2))=1+g_K.$$ The elements  
$
1,t,\dots,t^{g_K} 
$ are $\F_q$-linearly independent in $\sL_1(g_K(P_1+P_2))$, and therefore form a basis. Consequently, 
$$
\sL_1(r(P_1+P_2)) \subset \sL_1(g_K(P_1+P_2)) \subset \F_q[t],
$$ for every $r=0,\dots,g_K$. It follows that 
$$\sL_1(rP_i) \subset \sL_1(r(P_1+P_2)) \subset \F_q[t],\quad i=1,2,
$$ for every $r=0,\dots,g_K$. Every nonconstant polynomial has the same pole order at both $P_1$ and $P_2$, so
$$
\sL_1( rP_1)=\sL_1(r P_2)=\F_q,
$$ for every $r=0,\dots,g_K$. Therefore, by \cite[Theorem 1.6.8]{Sti09}, 
$$
W(P_i)=\N \setminus\{1,2,\dots,g_K\},\quad i=1,2.
$$
\end{proof}

We next count points with prescribed pole order. Let $K/k_0$ be a quadratic extension with $q_K=q$. If there is a unique place $P$ of $K$ above $\infty$, for $m \geq 0$, define
$$
\sL_n'(m \cdot P):= \{  \alpha \in K^n \mid (\alpha)_\infty=m \cdot P\}=\{ \alpha \in \cO_K^n \mid \max\{0,-\ord_P(\alpha)\}=m \}.
$$ If $\infty$ splits into two places $P_1$ and $P_2$ of $K$, for $a_1,a_2 \geq 0$, define
$$
\sL_n'(a_1P_1+a_2P_2)=\{ \alpha \in K^n \mid (\alpha)_\infty=a_1 P_1+a_2 P_2\}.
$$

\begin{lem}\label{clm:totram1} 
Suppose that $\infty$ is ramified in $K$, and let $P$ be the unique place above $\infty$. Then
\begin{equation*}
\#\sL_n'(m \cdot P)=
\begin{cases}
q^n &\textrm{if } m=0,\\
(q^n-1)q^{\frac{mn}{2}}&\textrm{if }0<m \leq 2g_K\textrm{ and }m\textrm{ is even},\\
(q^n-1)q^{n(m-g_K)}&\textrm{if }m \geq 2g_K+1,\\
0&\textrm{otherwise}.
\end{cases}
\end{equation*}
\end{lem}

\begin{proof} The case $m=0$ is immediate. For $m \geq 1$, we have
\begin{equation}\label{eq:diffL'}
\#\sL_n'(m \cdot P)=\#\sL_n(m\cdot P)-\#\sL_n((m-1)\cdot P).
\end{equation} 
For $0 \leq r \leq g_K$, the description of $W(P)$ in Lemma \ref{lem:WeieGap}, along with Clifford's theorem implies
$$
\ell_1(2r\cdot P)=r+1, \quad \ell_1((2r-1)\cdot P)=r\,\,(r\geq 1).
$$
Hence
$$
\#\sL_n(2r\cdot P)=q^{n(r+1)}\quad\textrm{and}\quad\#\sL_n((2r-1)\cdot P)=q^{nr}.
$$
Equation \eqref{eq:diffL'} gives the stated formula for $m \leq 2g_K$. If $m \geq 2g_K+1$, Equation \eqref{eq:diffL'} and Riemann--Roch applied to $m P$ and $(m-1)P$ yields
$$
\#\sL'_n(mP)=(q^n-1)q^{n(m-g_K)}.
$$
\end{proof}

\begin{lem}\label{clm:inert1} Suppose that $\infty$ is inert in $K$, and let $P$ be the unique place above $\infty$. Then
$$
\#\sL_n'(m \cdot P)=
\begin{cases}
q^n &\textrm{if } m=0,\\
(q^n-1)q^{nm} &\textrm{if }0< m \leq g_K,\\
(q^{2n}-1)q^{n(2m-g_K-1)}&\textrm{if } m \geq g_K+1.
\end{cases}
$$ 
\end{lem}

\begin{proof} Since $\deg(P)=2$, Lemma \ref{lem:bounddima} and Lemma \ref{lem:WeieGap}, imply that, for each $r \geq 1$,
\begin{equation}\label{eq:elelr}
\ell_1(rP)-\ell_1((r-1)P):=\varepsilon_r \in \{1,2\}.
\end{equation}
Riemann--Roch gives $
\ell_1(g_KP)=1+g_K$, whereas $\ell_1(0)=1$. Thus
\begin{align*}
&\sum_{r=1}^{g_K} \varepsilon_r=g_K,
\text{ and thus }\\
&\varepsilon_r=1 ,( r=1,\ldots, g_K).
\end{align*}
Therefore using \eqref{rieroch_n}, we have
$$
\#\sL_n(rP)=q^{n(r+1)} \quad(r=0,\dots,g_K).
$$
The first two cases now follow by taking successive differences, with $m=0$ treated separately.

For $m \geq g_K+1$, Riemann--Roch gives
$$
\#\sL_n'(m \cdot P)= \#\sL_n(m \cdot P)-\#\sL_n((m-1) \cdot P)=(q^{2n}-1)q^{n(2m-g_K-1)}.
$$
\end{proof}

\begin{lem}\label{clm:riemrochdimeP1P2} Suppose that $\infty$ splits into two places $P_1$ and $P_2$ of $K$. Then
$$
\#\sL_n(a_1P_1+a_2P_2)=
\begin{cases}
q^{n(a_1+a_2-g_K+1)}&\textrm{if }a_1 \geq g_K\textrm{ or }a_2 \geq g_K,\\
q^{n(\min\{a_1,a_2\}+1)}&\textrm{if }a_1,a_2\leq g_K,
\end{cases}
$$ and

$$\#\sL_n'(a_1P_1+a_2P_2)=$$
\begin{alignat*}1
\begin{cases}
q^{n(a_1+a_2-g_K-1)}(q^n-1)^2&\textrm{if }a_1 \geq g_K+1\textrm{ or }a_2 \geq g_K+1,a_1\cdot a_2\neq 0,\\
q^{n(a_1+a_2-g_K)}(q^n-1)&\textrm{if }a_1 \geq g_K+1\textrm{ or }a_2 \geq g_K+1, a_1\cdot a_2= 0,\\
q^{n a_1}(q^n-1)&\textrm{if }0<a_2=a_1\leq g_K,\\
q^n &\textrm{if } a_1=a_2=0,\\
0&\textrm{if }a_1\neq a_2\textrm{ and }a_1,a_2\leq g_K.
\end{cases}
\end{alignat*}
\end{lem}

\begin{proof}
By Lemma \ref{lem:WeieGap}(2), and the fact that $\deg P_1=\deg P_2=1$,
\begin{align*}
\ell_1(a_1P_1+a_2P_2)&=\max\{a_1,g_K\}+\max\{a_2,g_K\}-2g_K+\ell_1(\min\{a_1,g_K\}P_1+\min\{a_2,g_K\}P_2)\\
&=\max\{a_1,g_K\}+\max\{a_2,g_K\}-2g_K+\min\{a_1,a_2,g_K\}+1.
\end{align*}
This gives the first formula.

If $a_1,a_2>0$, inclusion-exclusion gives
\begin{align*}
\#\sL_n'(a_1P_1+a_2P_2)&=\#\sL_n(a_1P_1+a_2P_2)+\#\sL_n((a_1-1)P_1+(a_2-1)P_2)\\
&-\#\sL_n((a_1-1)P_1+a_2P_2)-\#\sL(a_1P_1+(a_2-1)P_2).
\end{align*}
If $a_i >0$ and the other coefficient is zero, then
\begin{align*}
    \#\sL_n'(a_iP_i)&=\#\sL_n(a_iP_i)-\#\sL_n((a_i-1)P_i).
\end{align*}
Substitution of the first formula yields the result.
\end{proof}

\section{Proof of Theorem \ref{thm:poinbounheiK}}\label{sec:poibouhei}
We prove Theorem \ref{thm:poinbounheiK} over an arbitrary function field. We use several results from \cite{BDW26}. Recall that $\cM^\infty_K$ is the set of places of $K$ above $\infty$, 
$$s_K=\#\cM^\infty_K, \quad \ell=[K:\F_q(t)],$$ and that $q_K=q$. Fix $n \geq 1$ and set  $\cA_n:=\cO_K^n$.  For $m \in U_K$,\footnote{By Remark \ref{rmk:heihpropre}, the values of $h_K$ on $\cO_K^n$ lie in
$
 U_K=\sum_{v \in \cM^\infty_K}\N \cdot \deg(v)
$.} let 
$$
\cA_{n,m}:=\left\lbrace\gamma \in \cA_n \mid h(\gamma)=\frac{m}{\ell}\right\rbrace=\{ \gamma \in \cA_n \mid h_K(\gamma)=m \}.
$$ 
Define
$$
\mathbb{B}_m:=\left\lbrace \underline{b}=(b_v)_{v \in \cM_K^\infty} \in \N^{s_K} \bigg\vert \sum_{v \in M^\infty_K} b_v \deg(v)=m \right\rbrace.
$$
For $\underline{b} \in \N^{s_K}$, set
$$
\cA_{n,\underline{b}}=\left\lbrace \gamma \in \cA_n \bigg\vert (\gamma)_{\infty}=\sum_{v \in M^\infty_K} b_v \cdot v\right\rbrace.
$$
For each
$\gamma \in \cA_{n,m}$, there is a unique 
$\underline{b} \in \mathbb{B}_m$ such that $\gamma \in \cA_{n,\underline{b}}$, namely
$$
b_v=\max\{0,-\ord_v(\gamma)\}.
$$
Thus
\begin{equation}
\cA_{n,m}=\bigsqcup_{\underline{b} \in \mathbb{B}_m}\cA_{n,\underline{b}},
\end{equation}
and hence
\begin{equation}\label{eq:33}
\#\cA_{n,m}=\sum_{\underline{b} \in \mathbb{B}_m} \#\cA_{n,\underline{b}}.
\end{equation}

\subsubsection*{Asymptotics of $\mathbb{B}_m$} For $\underline{b} \in \mathbb{B}_m$, set 
$$
\cM^\infty_{\underline{b}}:=\{ v \in \cM^\infty_K \mid b_v>0\}.
$$ Partition $\mathbb{B}_{m}$ as

\begin{align}\label{eq:defb00}
\mathbb{B}_{m}^{(0)}=\{ \underline{b} \in \mathbb{B}_m  \mid \cM^\infty_{\underline{b}}=\cM^\infty_K\}\quad\textrm{and}\quad \mathbb{B}_m^{(1)}=\mathbb{B}_m \setminus \mathbb{B}_m^{(0)}. 
\end{align}
Write 
$$
\cM^\infty_K=\{ v_1,\dots,v_{s_K}\} \quad \textrm{and} \quad \rho_K=\gcd(\deg(v_1),\dots,\deg(v_{s_K})).
$$ The following estimates are proved in \cite[\S 3.2]{BDW26}.
\begin{lem}\label{lem:dedea}
We have
$$
\# \mathbb{B}_m=\frac{\rho_K}{(s_K-1)! \prod_{v \in \cM^\infty_K}\deg(v)} \cdot m^{s_K-1}+ O_{s_K}\left( 1 + m^{s_K-2} \right),
$$
$$
\# \mathbb{B}_m^{(1)}=O_{s_K}(1+m^{s_K-2}),
$$
and
$$
\#\mathbb{B}_m^{(0)}=\frac{\rho_K}{(s_K-1)! \prod_{v \in \cM^\infty_K} \deg(v)}\cdot m^{s_K-1}+ O_{s_K} (1+m^{s_K-2}).
$$
\end{lem}

\subsubsection*{Asymptotics of  $\cA_{n,\underline{b}}$} 
Recall that 
$$
D_K=\sum_{v \in \cM^\infty_K}\deg(v).
$$ For $\underline{b} \in \mathbb{B}_m$, set 
$$
Q_{\underline{b}}=\sum_{v \in \cM^\infty_K}b_v v.
$$

\begin{lem}\label{lem:Bb} \label{cor:22as} Let $\underline{b} \in \mathbb{B}^{(0)}_m$. Then 
$$
\#\cA_{n,\underline{b}}=\left( \prod_{v \in \cM^\infty_K} \left(1 -\frac{1}{q^{n\deg(v)}}\right) \right) \cdot q^{mn-ng_K+n}+
\begin{cases}
0&\textrm{if }m\geq 2 g_K-1+D_K,\\
O_{s_K}\left(q^{\frac{m}{2}n+n}\right)&\textrm{otherwise}.
\end{cases}
$$
\end{lem}

\begin{proof} Since $\cM^\infty_K=\cM^\infty_{\underline{b}}$, we have
$$
\cA_{n,\underline b}
=
\left\{ \alpha \in K^n \mid
\ord_v(\alpha) = -b_v \quad \text{for all } v \in \cM^\infty_K \right\}.
$$ 
For $I \subset \cM^\infty_K$, let $e_I=\sum_{v \in I} v \in \cM_K$. Then
$$
\cA_{n,\underline{b}}=\mathscr{L}_n(Q_{\underline{b}}) \setminus \left( \bigcup_{v \in \cM^\infty_K } \mathscr{L}_n(Q_{\underline{b}}-e_{\{v\}})\right).
$$
Inclusion-exclusion gives 
$$
\# \cA_{n,\underline{b}}=\sum_{I \subset \cM^\infty_K} (-1)^{\# I} \# \mathscr{L}_n(Q_{\underline{b}}-e_{I}). 
$$  
Moreover, for $I \subset \cM^\infty_K$, we have
\begin{equation*}\label{eq:Bund}
\deg(Q_{\underline{b}}-e_{I})=m- \sum_{v \in I} \deg(v). 
\end{equation*}
\noindent
$\bullet$ Suppose that $\deg(Q_{\underline{b}}-e_{I}) < 2 g_K-1$. Then Clifford's theorem gives
$$
\#\mathscr{L}_n(Q_{\underline{b}}-e_I) \leq  q^{\frac{1}{2}n\left(m-\sum_{v \in I}\deg(v)\right)+n} \leq  q^{\frac{m}{2}n+n}.
$$ Moreover $$
q^{n(m-\sum_{v \in I} \deg(v)-g_K+1)} \leq q^{\frac{n}{2}(m-\sum_{v \in I}\deg(v))} \leq q^{\frac{nm}{2}}.
$$
\vskip 1mm
\noindent
$\bullet$ Suppose that $\deg(Q_{\underline{b}}-e_{I}) \geq 2 g_K-1$. Then the Riemann--Roch theorem implies
$$
\#\mathscr{L}_n(Q_{\underline{b}}-e_I)=q^{n(\deg(Q_{\underline{b}}-e_I)-g_K+1)}=q^{m n-n\sum_{v \in I} \deg(v)-n g_K+n}.
$$
Combining the above estimates, we obtain
\begin{align*}
\#\cA_{n,\underline{b}}&=\sum_{I \subset \cM^\infty_K} (-1)^{\#I} \# \mathscr{L}_n(Q_{\underline{b}}-e_I)\\
&=\sum_{I \subset \cM^\infty_K} \left( (-1)^{\#I} \left[ q^{-n \sum_{v \in I} \deg(v)} \cdot q^{mn-ng_K+n} \right]\right)+ \begin{cases}
0&\textrm{if }m\geq 2 g_K-1+D_K,\\
O_{s_K}\left(q^{\frac{m}{2}n+n}\right)&\textrm{otherwise},
\end{cases}\\
&=\left( \prod_{v \in \cM^\infty_K} \left( 1 - \frac{1}{q^{n \deg(v)}} \right) \right) \cdot q^{mn-n g_K+n}+
\begin{cases}
0&\textrm{if }m\geq 2 g_K-1+D_K,\\
O_{s_K}\left(q^{\frac{m}{2}n+n}\right)&\textrm{otherwise}.
\end{cases}
\end{align*}
\end{proof}

\begin{lem}\label{lem:gengen334}
Let $\underline{b} \in \mathbb{B}_m$. Then 
$$
\#\cA_{n,\underline{b}} \leq
\begin{cases}
q^{mn+n-ng_K}&\textrm{if }m \geq 2g_K-1,\\
q^{\frac{m}{2}n+n}&\textrm{otherwise}.
\end{cases}
$$
\end{lem}

\begin{proof}
Since $\cA_{n,\underline{b}} \subset \mathscr{L}_n(Q_{\underline{b}})$ and $\deg(Q_{\underline{b}})=m$, the desired estimate follows immediately from the Riemann--Roch theorem and Clifford's theorem.
\end{proof} 

\subsubsection{Completion of the proof of Theorem \ref{thm:poinbounheiK}}
By \eqref{eq:33} and \eqref{eq:defb00},
$$
N_{\cO_K}(n,m)=\sum_{\underline{b} \in \mathbb{B}_m} \#\cA_{n,\underline{b}}=\sum_{\underline{b} \in \mathbb{B}^{(0)}_m} \#\cA_{n,\underline{b}}+\sum_{\underline{b} \in \mathbb{B}^{(1)}_m} \#\cA_{n,\underline{b}}.
$$

Suppose first that $s_K \geq 2$. Applying Lemma \ref{lem:Bb} to the first sum and Lemma \ref{lem:gengen334} to the second one, and then using Lemma \ref{lem:dedea}, gives
$$
N_{\cO_K}(n,m)=S_{\cO_K}(n) m^{s_K-1}q^{mn}+O_{s_K}(E(m)),
$$
where $E(m)$ exactly as stated in Theorem \ref{thm:poinbounheiK}.

Now suppose that $s_K=1$ and write $\cM_K^\infty=\{w\}$. Since
 $\mathbb{B}_m=\{b= m /\deg(w)\}$ and $m>0$, we have $\cM^\infty_{b}=\cM^\infty_K$. Therefore, the desired formula follows directly from Lemma \ref{cor:22as}.

\section{Proof of Theorem \ref{thm:main1}}\label{ss:prfmain1}
We begin by proving Corollary \ref{cor:ofdegr}, an analogue of \cite[Corollary 2]{TW13}. Recall that 
$$
d=[k:\F_q(t)] \quad \textrm{and}\quad q_k=q.
$$ For every extension $K/k$ of degree $e$ with $q_K=q$, recall
$$
\delta_{n+1}(K/k)=\min\{h(P)\mid P \in \mathbb{P}^n(K)\,\,\textrm{and}\,\, k(P)=K\}.
$$
Throughout this section, we use the natural embedding 
$$
(\alpha_1:\dots:\alpha_n) \in \overline{k}^n \mapsto (1:\alpha_1:\dots:\alpha_n) \in \mathbb{P}^n(\overline{k}).
$$
\begin{proof}[Proof of Corollary \ref{cor:ofdegr}]
For $m \in U_K$,
\begin{equation}\label{eq:NOkl}
N_{\cO_k}(n,K,m)
= N_{\cO_K}(n,m)
-\sum_{\substack{e' \mid e \\ e' < e}}
\ \sum_{\substack{k \subset L \subset K \\ [L:k]=e'}}  N_{\cO_k}\left(n,L, \frac{me'}{e}\right).
\end{equation}
\noindent
For each $L$ occurring in \eqref{eq:NOkl}, embed $L^n$ into $\mathbb{P}^n(L)$ and apply \cite[Corollary 1]{TW13}. Since $e' \leq e/2$, 
$$
N_{\cO_k}(n,L, me'/e)\ll_{n,k,e} q^{\frac{m(n+1)}{2}}.
$$ 
By \cite[Lemma 6]{TW13}\footnote{which bounds the number of intermediate fields.}, it follows that
$$
\sum_{\substack{e' \mid e \\ e' < e}}
\ \sum_{\substack{k \subset L \subset K \\ [L:k]=e'}} N_{\cO_k}(n, L,  me'/e)\ll_{n,k,e} q^{\frac{m(n+1)}{2}}.
$$
Combining this estimate with Theorem \ref{thm:poinbounheiK}, and using
$$
m^{s_K-1}q^{\frac{nm}{2}+n}=O_{k,n,e}(q^{\frac{m(n+1)}{2}}),
$$
gives the corollary.
\end{proof}

Towards the proof of Theorem \ref{thm:main1}, we now estimate the contribution from extensions of large genus, which corresponds to the error term appearing in $N_{\cO_k}(n,e,m)$. For a finite extension $K/k$ with $q_K=q$, let $J_K$ be the number of divisor classes of degree $0$, let $\zeta_K$ be the zeta function of $K$, and define the \emph{Schanuel constant}
$$
S_K(n)=\frac{J_K}{(q-1)\zeta_K(n) q^{n(g_K-1)}}.
$$

\begin{lem}\label{lem:boundmleq2g} Assume $3 \leq e < n$. Then
$$
\sum_{\substack{[K:k]=e\\q_K=q \\ m \leq 2g_K-2}} N_{\cO_k}(n,K,m)\ll_{n,k,e} 
\begin{cases}
mq^{mn}&\textrm{if }n=e+1,e \geq 5, \\
q^{mn} &\textrm{if }n > e+1.
\end{cases}
$$
Further, when $(n,e)\in \{(4,3),(5,4)\}$, we have for $\varepsilon>0$
$$
\sum_{\substack{[K:k]=e\\q_K=q \\ m \leq 2g_K-2}} N_{\cO_k}(e+1,K,m)\ll_{n,k,e,\varepsilon} 
\begin{cases}
q^{4m+m\varepsilon} &\textrm{if }e=3.\\
q^{5m+m\varepsilon+\frac{m}{3}}&\textrm{if }e= 4.
\end{cases}
$$
\end{lem}

\begin{proof} 

Let $\varepsilon>0$. As before, embed $K^n$ into $\mathbb{P}^{n}(K)$. By \cite[Theorem 4.1 (1) \& (2)]{Thu25} and \cite[Theorem 5.1]{Thu25}, we have
$$
\sum_{\substack{[K:k]=e\\q_K=q \\ m \leq 2g_K-2}} N_{\cO_k}(n,K,m) \ll_{n,k,e} \sum_{\substack{[K:k]=e\\q_K=q \\ m \leq 2g_K-2}} \frac{J_K}{q^{(n+1)g_K}}q^{(n+1)m}+ \begin{cases}
mq^{m(e+1)} &\textrm{if }n=e+1,\\
q^{mn} &\textrm{if } n>e+1.
\end{cases}
$$
On the other hand, by \cite[Lemma 6.1]{Thu25},
$$
\sum_{\substack{[K:k]=e\\q_K=q \\ m \leq 2g_K-2}} \frac{J_K}{q^{(n+1)g_K}}q^{(n+1)m} \ll_{n,k,e,\varepsilon}
\begin{cases}
q^{\frac{m}{2}(n+4+ \varepsilon)}&\textrm{if }e=3,\\
q^{\frac{m}{2}\left(n+\frac{e}{2}+3+\varepsilon\right)}&\textrm{if }e>3\textrm{ and }p\nmid e,\\
q^{\frac{m}{2} \left(n+\frac{e}{2}+\frac{1}{e-1}+3+\varepsilon \right)}&\textrm{if }e>3\textrm{ and }p\mid e,
\end{cases}
$$
where $p$ is the characteristic of $K$. \\
The second part of the Lemma follows from the above estimates. 

For the first part, we take $\varepsilon=\frac{1}{4}$. Observing that $e\ge 3$, we obtain
\begin{align}\label{eq: sumhighgen}
\sum_{\substack{[K:k]=e\\q_K=q \\ m \leq 2g_K-2}} N_{\cO_k}(n,K,m)\ll_{n,k,e} 
\begin{cases}
mq^{mn}&\textrm{if }n=e+1,e \geq 5, \\
q^{mn} &\textrm{if }n > e+1.
\end{cases}
\end{align}

\end{proof}

The quadratic analogue follows from \cite{Thu15}.
\begin{lem}\label{lem:boundmleq2gz} Assume that $e=2$ and $n \geq 3$. Then
\begin{align*}
\sum_{\substack{[K:k]=2\\ q_K=q \\ m \leq 2g_K-2}}N_{\cO_k}(n,K,m)\ll_{n,k} q^{(n/2+2)m}.
\end{align*}
\end{lem}

\begin{proof}
By \cite[Lemma 7]{TW13}, there is a constant $c_k>0$, depending only on $k$, such that $N_{\cO_K}(n,m)=0$ whenever $K/k$ is a quadratic, $q_K=q$ and $g_K>m+c_k$. By \cite[Theorem 1]{Thu15}, 
\begin{align*}
N_{\cO_k}(n,K,m)\ll_{n,k} 
S_K(n+1)q^{(n+1)m}+O_{n,k}\left(q^{n(m-g_K)}\left(J_K+q^{m/2}\sqrt{J_K \log_q^+(J_K)}\right)\right).
\end{align*}
We sum the three terms separately over
$$
\frac{m+2}{2} \leq g_K \leq m+c_k.
$$
First, \cite[Proposition 3]{Thu15} gives, for $n \geq 3$,
$$
\sum_{\substack{[K:k]=2\\ q_K=q\\\frac{m+2}{2} \leq g_K \leq m+c_k}} S_K(n+1)\ll_{n,k} \sum_{\substack{[K:k]=2\\ q_K=q\\\frac{m+2}{2} \leq g_K \leq m+c_k}} J_K q^{-(n+1)g_K}\ll_{n,k} \sum_{i\geq \frac{m+2}{2}} q^{-(n-2)i}\ll_{n,k} q^{-(n-2)\frac{m}{2}}.
$$ 
Next, by \cite[Equation (28)]{Thu15},
$$
\sum_{\substack{[K:k]=2\\ q_K=q\\\frac{m+2}{2} \leq g_K \leq m+c_k}} q^{n(m-g_K)}J_K \ll_{n,k}
\begin{cases}
m q^{nm} &\textrm{if } n=3,\\
q^{\left(\frac{n+3}{2}\right)m} &\textrm{if }n> 3.
\end{cases}
$$
Finally, \cite[Equation (29)]{Thu15} gives
\begin{align*}
\sum_{\substack{[K:k]=2\\ q_K=q\\\frac{m+2}{2} \leq g_K \leq m+c_k}} q^{n(m-g_K)} q^{m/2} \sqrt{J_K \log_q^+(J_K)} \ll_{n,k}\sqrt{m} q^{\left(n+\frac{1}{2}\right)m} q^{- \left(n-\frac{5}{2} \right)\frac{m}{2}}\ll_{n,k} \sqrt{m} q^{\left(\frac{n}{2} + \frac{7}{4} \right)m}.
\end{align*}
Combining these estimates, we obtain
\begin{align*}
\sum_{\substack{[K:k]=2\\ q_K=q \\ m \leq 2g_K-2}}N_{\cO_k}(n,K,m)&\ll_{n,k} q^{(n/2+2)m}+ \sqrt{m} q^{\left(\frac{n}{2} + \frac{7}{4} \right)m}+\begin{cases}
m q^{nm} &\textrm{if } n=3,\\
q^{\left(\frac{n+3}{2}\right)m} &\textrm{if }n> 3.
\end{cases}\\
&\ll_{n,k} q^{(n/2+2)m}.
\end{align*}

\end{proof}

\subsection{Completion of the proof of Theorem \ref{thm:main1}}
For every degree-$e$ extension $K/k$ with $q_K=q$, recall that
$$
S_{\cO_K}(n)= \kappa_n(K) q^{-ng_K+n}.
$$
Moreover
$$
S_{\cO_k}(n,e) = \sum_{K} S_{\mathcal{O}_K}(n),
$$
where the sum is over the degree-$e$ extensions $K/k$ with $q_K=q$ in which every place above $\infty$ splits completely. We also have
\begin{equation}
N_{\cO_k}(n,e,m)
= \sum_{\substack{[K:k]=e \\ q_K=q\\ \delta_{n+1}(K/k) \leq m/(de)}} N_{\cO_k}(n,K,m).
\end{equation}
For $s_k \leq i \leq e s_k$, let $\sC_{k,i}$ be the set of degree-$e$ extensions $K/k$ with $q_K=q$ for which $\#\cM_K^\infty=i$, and define
$$
N^{(i)}_{\cO_k}(n,e,m)
= \sum_{\substack{K \in \sC_{k,i}\\ \delta_{n+1}(K/k) \leq m/(de)}} N_{\cO_k}(n,K,m).
$$
Then,
\begin{equation}\label{eq:proof:thmmain1}
N_{\cO_k}(n,e,m)=\sum_{i=s_k}^{es_k} N^{(i)}_{\cO_k}(n,e,m).
\end{equation}

\begin{claim}\label{clm:leveli}
Let $s_k \leq i \leq e s_k$ and let $\varepsilon>0$. 
\begin{enumerate}
    \item We have
\begin{align*}
N^{(i)}_{\cO_k}(n,e,m)& \leq m^{i-1}q^{mn} \sum_{\substack{ K \in \sC_{k,i} \\ \delta_{n+1}(K/k) \leq m/(de) \\ 2g_K-1 \leq m}}S_{\cO_K}(n)+O_{k,n,e} \left(m^{i-2}\, q^{mn}+q^{\frac{m(3n+4)}{4}}\right)\\
&+ \begin{cases}
O_{k,n,e}(mq^{mn})&\textrm{if }n=e+1\geq 6,\\
O_{k,n,e}(q^{mn}) &\textrm{if }n > e+1\geq 4,\\
O_{k,n,e}\left(q^{\left(\frac{n}{2}+2\right)m}\right) &\textrm{if }n>e+1,e=2.
\end{cases}
\end{align*}
Moreover, equality holds if $i=es_k$.

\item  Additionally, 
\begin{align*}
    N^{(i)}_{\cO_k}(e+1,e,m)&=\begin{cases}
        O_{k,e}(q^\frac{7m}{2}) &\textrm{ if } e=2\\
        O_{k,e,\varepsilon}(q^{4m+m\varepsilon}) &\textrm{if }e=3.\\
        O_{k,e,\varepsilon}(q^{5m+m\varepsilon+\frac{m}{3}})&\textrm{if }e= 4.
    \end{cases}
\end{align*}
\end{enumerate}
\end{claim}

\begin{proof}[Proof of Claim \ref{clm:leveli}]

We have
\begin{align}\label{eq:sum_decomp}
\nonumber N^{(i)}_{\cO_k}(n,e,m)
&= \sum_{\substack{ K \in \sC_{k,i} \\ \delta_{n+1}(K/k) \leq m/(de)}}
   N_{\cO_k}(n,K,m) \\
&= \sum_{\substack{ K \in \sC_{k,i} \\ \delta_{n+1}(K/k) \leq m/(de) \\ 2g_K-1 \leq m}}
   N_{\cO_k}(n,K,m)
 + \sum_{\substack{ K \in \sC_{k,i} \\ \delta_{n+1}(K/k) \leq m/(de) \\ m \leq 2g_K-2}}
   N_{\cO_k}(n,K,m).
\end{align}

For the first term, by Corollary \ref{cor:ofdegr}, we obtain \footnote{In fact, for $K \in \sC_{k,es_k}$, we have $U_K=U_k$.}
\begin{align}\label{eq: first_main_sum}
\nonumber\sum_{\substack{ K \in \sC_{k,i} \\ \delta_{n+1}(K/k) \leq m/(de) \\ 2g_K-1 \leq m}}
   N_{\cO_k}(n,K,m)&\leq m^{i-1}q^{mn} \sum_{\substack{ K \in \sC_{k,i} \\ \delta_{n+1}(K/k) \leq m/(de) \\ 2g_K-1 \leq m}}S_{\cO_K}(n) \quad\textrm{(with equality when $i=es_k$)}\\
   \nonumber&+O_{k,n,e} \left(\sum_{\substack{ K \in \sC_{k,i} \\ \delta_{n+1}(K/k) \leq m/(de) \\ 2g_K-1 \leq m}}\left(m^{i-2}\, q^{mn-ng_K}+q^{\frac{m(n+1)}{2}}\right)\right)\\
   \nonumber&=m^{i-1}q^{mn} \sum_{\substack{ K \in \sC_{k,i} \\ \delta_{n+1}(K/k) \leq m/(de) \\ 2g_K-1 \leq m}}S_{\cO_K}(n)+O_{k,n,e} \left(m^{i-2}\, q^{mn}\right)\\
   &+
   O_{k,n,e}\left(q^{\frac{m(3n+4)}{4}}\right).
\end{align}
The last equality in (\ref{eq: first_main_sum}) follows from the following. By \cite[Theorem]{Thu13} and \cite[Proposition 3]{Thu15}, we have
\begin{equation}\label{eq:bounddd}
\sum_{\substack{[K:k]=e\\ q_K=q\\g_K=j}}1 \ll_{k,e}
\begin{cases}
q^{j\left( \frac{e+2}{2}+\frac{1}{e-1} \right)} &\textrm{if }p\mid e,e>2,\\
q^{j\frac{e+2}{2}}&\textrm{if } p\nmid e, e>2,\\
q^{2j}&\textrm{ if } e=2,
\end{cases}
\end{equation}
where $p$ denotes the characteristic of $\F_q$.
It follows that, for $n \geq 3$ and $n  \geq e+1$, 
$$
\sum_{\substack{ K \in \sC_{k,i} \\ \delta_{n+1}(K/k) \leq m/(de) \\ 2g_K-1 \leq m}} q^{-n g_K} \ll_{n,k,e} 1,
$$
and 
$$
\sum_{\substack{ K \in \sC_{k,i} \\ \delta_{n+1}(K/k) \leq m/(de) \\ 2g_K-1 \leq m}} q^{\frac{m(n+1)}{2}} \ll_{n,k,e} 
\begin{cases}
q^{\frac{m(3n+4)}{4}} &\textrm{if }e\geq 3,\\
q^{\frac{m(n+3)}{2}}&\textrm{if }e=2.
\end{cases}
$$
Combining these estimates,  observing that $n+3\le \frac{3n+4}{2}$ for all $n\ge 2$, yields (\ref{eq: first_main_sum}).

For the second term in (\ref{eq:sum_decomp}), for $\varepsilon>0$, by Lemma \ref{lem:boundmleq2g} and Lemma \ref{lem:boundmleq2gz}, we have
\begin{equation*}
\sum_{\substack{ K \in \sC_{k,i}  \\ \delta_{n+1}(K/k) \leq m/(de) \\ m \leq 2g_K-2}}
   N_{\cO_k}(n,K,m) \ll_{n,k,e} \begin{cases}
mq^{mn}&\textrm{if }n=e+1,e\geq 5,\\
q^{mn}&\textrm{if }n>e+1,e\geq 3,\\
q^{(n/2+2)m}&\textrm{if }n\ge e+1,e=2,
\end{cases}
\end{equation*}

and,
$$
\sum_{\substack{[K:k]=e\\q_K=q \\ m \leq 2g_K-2}} N_{\cO_k}(e+1,K,m)\ll_{n,k,e,\varepsilon} 
\begin{cases}
q^{4m+m\varepsilon} &\textrm{if }e=3.\\
q^{5m+m\varepsilon+\frac{m}{3}}&\textrm{if }e= 4.
\end{cases}
$$

These estimates, together with (\ref{eq: first_main_sum}) prove the claim.
\end{proof}

\begin{claim}\label{clm:leadterm}
We have
$$
 \sum_{\substack{K \in \sC_{k,s_ke}\\\delta_{n+1}(K/k) > m/(de) \textrm{ or } 2g_K-1 > m}}\kappa_n(K) q^{-ng_K}=
\begin{cases}
O_{k,n,e}(q^{\left(-n+\frac{e+2}{2}+\frac{1}{e-1}\right)\frac{m}{2}})&\textrm{if }e\geq 3,\\
 O_{k,n,e}(q^{(-n+2)m/2}) &\textrm{if }e=2.
\end{cases}
$$
\end{claim}

\begin{proof}[Proof of Claim \ref{clm:leadterm}] Let $c_2(k,e)$ be the constant from \cite[Lemma 7]{TW13}. For every degree-$e$ extension $K/k$ with $q_K=q$,
$$
\delta_{n+1}(K/k) \leq \frac{g_K+c_2(k,e)}{de}.
$$
By Equation \eqref{eq:bounddd}, we obtain
\begin{align*}
\sum_{\substack{K \in \sC_{k,s_k e}\\\delta_{n+1}(K/k) > m/(de) \textrm{ or } 2g_K-1 > m}}\kappa_n(K) q^{-ng_K} &\le\sum_{\substack{K \in \sC_{k,s_k e}\\\delta_{n+1}(K/k) > m/(de)}}\kappa_n(K) q^{-ng_K}+ \sum_{\substack{K \in \sC_{k,s_k e}\\ 2g_K-1 > m}}\kappa_n(K) q^{-ng_K}  \\
& \leq \sum_{\substack{K \in \sC_{k,s_k e}\\ g_K \geq m-c_2(k,e)}}\kappa_n(K) q^{-ng_K}+ \sum_{\substack{K \in \sC_{k,s_k e}\\ 2g_K-1 > m}}\kappa_n(K) q^{-ng_K} \\
&\leq  
\begin{cases}
2\sum_{j \geq (m+1)/2} q^{-nj}q^{j\left( \frac{e+2}{2}+\frac{1}{e-1} \right)} &\textrm{if }e\geq 3,\\
2\sum_{j \geq (m+1)/2} q^{-nj}q^{2j} &\textrm{if }e=2.
\end{cases}
\\
&
\ll_{n,k,e} 
\begin{cases}
q^{\left(-n+\frac{e+2}{2}+\frac{1}{e-1}\right)\frac{m}{2}}&\textrm{if }e\geq 3,\\
q^{\frac{(-n+2)m}{2}}&\textrm{if }e=2.
\end{cases}
\end{align*}
This proves the claim.
\end{proof}
\noindent
By Claim \ref{clm:leadterm}, the series $S_{\cO_k}(n,e)$ converges for all $n\ge e+1\ge3$. Moreover, combining Claim \ref{clm:leveli} and Claim \ref{clm:leadterm}, we obtain the following:
 \begin{align*}
N_{\cO_k}(n,e,m)&=S_{\cO_k}(n,e) m^{s_ke-1}q^{mn}+O_{k,n,e} \left(m^{s_ke-1}q^{\frac{(n+3)m}{2}}+m^{s_ke-2}\, q^{mn}\right)\\ 
&+ \begin{cases}
O_{k,n,e}\left(mq^{mn}\right)&\textrm{if }n=e+1,e\geq 5,\\
O_{k,n,e}\left(q^{mn}\right)&\textrm{if }n>e+1,e\geq 3,\\
O_{k,n,e}\left(q^{(n/2+2)m}\right)&\textrm{if }n \geq 4,e=2.
\end{cases}
\end{align*}
\begin{align*}
N_{\cO_k}(e+1,e,m)&=
 \begin{cases}
O_{k,e}\left(q^{\frac{7m}{2}}\right)&\textrm{if }e=2,\\
O_{k,e, \varepsilon}\left(q^{4m+m\varepsilon}\right)&\textrm{if }e= 3,\\
O_{k,e,\varepsilon}\left(q^{5m+\frac{m}{3}+m\varepsilon}\right)&\textrm{if }e=4.
\end{cases}
\end{align*}

This proves Theorem \ref{thm:main1}.\\

\section{Proof of Theorem \ref{thm:main2}}\label{sec:smalcas}
In this section, we establish the estimates in Theorem \ref{thm:main2} . We assume throughout that $k=\F_q(t)$, where $q$ odd. The cases $n=2$ and $n=3$ use Weierstrass semigroups calculations from \S\ref{ssec:weismg}; the case $n=1$ is elementary. 

\subsection{Fixed quadratic extensions}
Throughout this subsection, $K$ denotes a quadratic extension of $k=\F_q(t)$, with $q_K=q$.

\begin{lem}\label{lem:inert_rami}
Let $n \geq 1$. If $\infty $ is inert or ramified in $K/k$, then 
 $$
 N_{\cO_k}(n,K,m)=0
 $$ for every $m \in U_K$ with $m\le 2g_K-2$. 
\end{lem}
\begin{proof} Let $P$ be the unique place of $K$ above $\infty$. Set 
$$
V_m=
\begin{cases}
\sL\left( \frac{m}{2} \cdot P \right)&\textrm{if }\infty \textrm{ is inert},\\
\sL\left( m \cdot P \right)&\textrm{if }\infty\textrm{ is ramified}. 
\end{cases}
$$
In the inert case, $m \in 2 \N$. In both cases, we have by Clifford's theorem, 
$$
\dim_{\F_q}(V_m) \leq 1+ \left\lfloor \frac{m}{2}\right\rfloor.
$$ The elements 
$$1,\dots, t^{^{\left\lfloor \frac{m}{2} \right\rfloor}},$$ are linearly independent in $V_m$, and hence form an $\F_q$-basis. Thus $V_m \subset \F_q[t]$. Every element $\alpha \in \cO_K^n$ of height $m/2$ has all its coordinates in $V_m$, and hence $\alpha \in \F_q[t]^n$, so $k(\alpha)=k$. Therefore, $N_{\cO_k}(n,K,m)=0$.  
\end{proof}

\begin{lem}\label{lem:Lspace_low height}\label{lem:split}
     Assume that $\infty $ splits into two places $P_1,P_2$ in $K$. Then, for $n \geq 1$ and $0\le b_1,b_2\le g_K$, 
 $$\sL_n(b_1P_1+b_2P_2)\subset \F_q[t]^n.$$ Consequently, $$
 N_{\cO_k}(n,K,m)=0
 $$ for all $m\le g_K$.
\end{lem}

\begin{proof}
It is enough to prove the first assertion for $n=1$. By Lemma \ref{clm:riemrochdimeP1P2}, 
$$
\ell_1(b_1P_1+b_2P_2)=1+\min\{b_1,b_2\}.
$$
Therefore $1,t,t^2,\dots t^{\min\{b_1,b_2\}}$  form an $\F_q$-basis of $\sL_1(b_1P_1+b_2P_2)$, which proves the inclusion. 
If $m \leq g_K$ and 
$$(\alpha)_\infty=b_1P_1+b_2P_2, \quad b_1+b_2=m,
$$ then $b_1,b_2\leq g_K$, so $\alpha \in  \F_q[t]^n$ and $k(\alpha)=k$. Thus $N_{\cO_k}(n,K,m)=0$.
\end{proof}

\begin{lem}\label{lem:n=2_fixed_K_count}
 Define
    $$
    E_{n,m}=
        \begin{cases}O(q^\frac{mn}{2})&\textrm{if $m$ is even},\\
        0 &\textrm{otherwise}.
        \end{cases}
    $$ 
    Where the implied constant depends only on $n$.\\
    Let $m>4$ with $m \in U_K$.  
    \begin{enumerate}
       \item Assume that $m-1\ge g_K\ge  1+\frac{m}{2}$. Then 
       $$N_{\cO_k}(n,K,m)=\begin{cases}
           2(m-g_K+\frac{1}{q^n-1})S_{\cO_K}(n)\cdot q^{nm}  &\textrm{ if $\infty$ splits in $K$,}\\
           0 &\textrm{ if $\infty$ is inert or ramified in $K$. }
       \end{cases}$$
        
        \item Assume that $g_K< 1+\frac{m}{2}$. If $\infty$ splits in K, then
        $$
        N_{\cO_k}(n,K,m)=\begin{cases}
      \left(m+\frac{q^n+1}{q^n-1}\right)  S_{\cO_K}(n)\cdot q^{nm}+E_{n,m} &\textrm{ if $m>2g_K$},\\
      \left(m+\frac{2}{q^n-1}\right)  S_{\cO_K}(n)\cdot q^{nm} &\textrm{ if $m=2g_K$},\\
      \left(m+\frac{3-q^n}{q^n-1}\right)  S_{\cO_K}(n)\cdot q^{nm} &\textrm{ if $m=2g_K-1$}.
        \end{cases}
        $$ 
    If $\infty$ is ramified in $K$, then
         $$
        N_{\cO_k}(n,K,m)= \begin{cases}
            S_{\cO_K}(n)\cdot q^{nm}+E_{n,m} &\textrm{if $m\ge 2g_K$},\\
       0 &\textrm{if $m=2g_K-1$}.
        \end{cases}
        $$
    If $\infty$ is inert in $K$, then
        $$
        N_{\cO_k}(n,K,m)= \begin{cases}
            S_{\cO_K}(n)\cdot q^{nm}+E_{n,m} &\textrm{if $m\ge 2g_K+2$},\\
        \frac{q^n}{q^n+1}S_{\cO_K}(n)\cdot q^{nm}+E_{n,m} &\textrm{if $m=2g_K$}.
        \end{cases}
        $$ 
    \item Assume $g_K \geq m$. Then $N_{\cO_k}(n,K,m)=0$.
    \end{enumerate}
\end{lem}
\begin{proof} Part (3) follows from Lemmas \ref{lem:inert_rami} and \ref{lem:split}. It remains to prove parts (1) and (2)
\vskip 1mm
\noindent
$\bullet$ Assume first that $\infty$ splits into two places $P_1$ and $P_2$ of $K$. For each $(b_1,b_2)\in \N^2$ satisfying $b_1+b_2=m$, let
 $$
 \sL_n'(b_1 P_1+b_2 P_2)=\{\alpha=(\alpha_1,\dots\alpha_n) \in \cO_K^n \mid (\alpha)_\infty=b_1P_1+b_2 P_2\}.
 $$
 Then,
 \begin{equation}\label{eq:heiLprim}
  N_{\cO_k}(n,K,m)=\sum_{b_1+b_2=m}\# (\sL'_n(b_1 P_1+b_2 P_2)\setminus \cO_k^n).
 \end{equation}  
 
 \begin{claim}\label{clm:slpslpslp} Suppose that $m-1 \geq g_K \geq 1+\frac{m}{2}$. Then
 $$
N_{\cO_k}(n,K,m)= 2\left( m-g_K +\frac{1}{q^n-1} \right) S_{\cO_K}(n)q^{nm}.
$$
 \end{claim}
 \begin{proof}[Proof of Claim \ref{clm:slpslpslp}]
The key point is to establish the following. For $(b_1,b_2) \in \N^2$ with $b_1+b_2=m$, we have 
\begin{equation}\label{eq:formform4233}
\#( \sL'_n(b_1 P_1+b_2 P_2)\setminus\cO_k^n)=
\begin{cases}
0 &\textrm{if }\max\{b_1,b_2\} \leq g_K,\\
S_{\cO_K}(n)q^{nm}&\textrm{if }0<b_1,b_2\textrm{ and }\max\{b_1,b_2\} \geq g_K+1,\\
S_{\cO_K}(n) \cdot \frac{q^n}{q^n-1} \cdot q^{nm}&\textrm{if }b_1=0 \textrm{ or }b_2=0.
\end{cases}
\end{equation}
The first case follows from Lemma~\ref{lem:Lspace_low height}. Suppose now we are in the second case. Since $g_K \geq 1+\frac{m}{2}$, we have $b_1\neq b_2$, and therefore $\#(\sL'_n(b_1P_1+b_2P_2)\setminus \cO_k^n)=\#\sL'_n(b_1P_1+b_2P_2)$. By the inclusion-exclusion principle and Lemma \ref{clm:riemrochdimeP1P2}, 
\begin{align*}
\#\sL'_n(b_1 P_1+b_2P_2)&=\#\sL_n(b_1P_1+b_2 P_2)-\#\sL_n((b_1-1)P_1+b_2 P_2)-\#\sL_n(b_1P_1+(b_2-1)P_2)\\
&+\#\sL_n((b_1-1)P_1+(b_2-1) P_2)\\
&=q^{n(m+1-g_K)}-q^{n(m-g_K)}-q^{n(m-g_K)}+q^{n(m-1-g_K)}\\
&\left(1 - \frac{1}{q^n} \right)^2 q^{-ng_K+n} q^{nm}\\
&= S_{\cO_K}(n)q^{nm}.
\end{align*}
Now suppose we are in the third case. Applying the inclusion-exclusion principle once again, we obtain 
$$
\#\sL_n'(b_1P_1+b_2 P_2)=S_{\cO_K}(n)\cdot \frac{q^n}{q^n-1}q^{nm}.
$$
This proves  \eqref{eq:formform4233}.

Finally, combining \eqref{eq:formform4233} with \eqref{eq:heiLprim}, we obtain
$$
N_{\cO_k}(n,K,m)=2(m-g_K-1) S_{\cO_K}(n) q^{nm}+ 2 S_{\cO_K}(n)\cdot \frac{q^n}{q^n-1}q^{nm}= 2\left( m-g_K +\frac{1}{q^n-1} \right) S_{\cO_K}(n)q^{nm}.
$$
Indeed, there are exactly $2(m-g_K-1)$ pairs $(b_1,b_2)$ satisfying $b_1+b_2=m$, $b_1,b_2 >0$, and $\max\{b_1,b_2\} \geq g_K+1$. This completes the proof of the claim.
\end{proof}

\begin{claim}\label{clm:clmspl2} Suppose that $g_K < 1 + \frac{m}{2}$. Then
$$
N_{\cO_k}(n,K,m)=\begin{cases}
      \left(m+\frac{q^n+1}{q^n-1}\right)  S_{\cO_K}(n)\cdot q^{nm}+E_{n,m} &\textrm{ if $m>2g_K$},\\
      \left(m+\frac{2}{q^n-1}\right)  S_{\cO_K}(n)\cdot q^{nm} &\textrm{ if $m=2g_K$},\\
      \left(m+\frac{3-q^n}{q^n-1}\right)  S_{\cO_K}(n)\cdot q^{nm} &\textrm{ if $m=2g_K-1$}.
        \end{cases}
$$ 

\end{claim}

\begin{proof}[Proof of Claim \ref{clm:clmspl2}]
Fix $(b_1,b_2)\in \N^2$ with $b_1+b_2=m$. By Lemma \ref{lem:Lspace_low height}, if $\max\{b_1,b_2\} \le g_K$, then $\sL_n'(b_1P_1+b_2P_2)\setminus\cO_k^n=\emptyset$; since $m \geq 2g_K-1$, this case occurs only when 
$$
(b_1,b_2 \in \{(g_K,g_K),(g_K,g_K-1),(g_K-1,g_K)\}.
$$

In the sequel, assume that $\max\{b_1,b_2\} \geq g_K+1$. As in the proof  of Claim \ref{clm:slpslpslp},
\begin{align}\label{eq:lpieee}
\#\sL_n'(b_1 P_1+ b_2 P_2)=\begin{cases}
         S_{\cO_K}(n)q^{nm} &\textrm{ if $b_1,b_2>0$},\\
    S_{\cO_K}(n)\cdot \frac{q^n}{q^n-1}q^{nm} &\textrm{ if $b_1=0$ or $b_2=0$}.
    \end{cases}
    \end{align}
Note that, when $b_1 \neq b_2$, 
$$
\sL_n'(b_1 P_1+b_2 P_2)=\sL_n'(b_1 P_1+b_2 P_2) \setminus \cO_k^n,
$$ and when $b_1=b_2=\frac{m}{2}$--which implies $m>2g_K$ and that $m$ is even--, we have  
$$\#(\sL_n'(b_1 P_1+b_2 P_2) \cap \cO_k^n)=O(q^\frac{mn}{2}).$$ 
Combining these observations, \eqref{eq:lpieee} and \eqref{eq:heiLprim}, the claim follows.
\end{proof}
\vskip 1mm
\noindent
$\bullet$ Assume next that $\infty$ ramifies, and let $P$ be the unique place of $K/k$ lying above $\infty$. Note that
$$
N_{\cO_k}(n,K,m)=\#(\sL_n'(mP) \setminus \cO_k^n).
$$
If $m-1 \geq g_K \geq 1 + \frac{m}{2}$, then Lemma \ref{lem:inert_rami} implies that
$$
N_{\cO_k}(n,K,m)=0.
$$
Suppose now that $g_K< 1 + \frac{m}{2}$. By Lemma \ref{lem:WeieGap}, $N_{\cO_k}(n,K,2g_K-1)=0$. Hence it remains to consider the case $m\ge 2g_K$. By Lemma \ref{clm:totram1}, 
$$
\sL_n'(mP)=
 \begin{cases}
      S_{\cO_K}(n)\cdot q^{nm} &\textrm{ if $m\ge 2g_K$},\\
        0 &\textrm{ if $m=2g_K-1$}.
 \end{cases}
$$
Moreover, any  $(\alpha)\in \cO_K^n$ with  $(\alpha)_\infty=mP$ satisfies $\alpha_1,\dots\alpha_n\in k$ if and only if $m$ is even and $\alpha_1,\dots\alpha_n\in \cO_k=\F_q[t]$ have degree at most $\frac{m}{2}$, with at least one of them having degree exactly $\frac{m}{2}$. There are $O(q^\frac{mn}{2})$ such $\alpha$ when $m$ is even, and none when $m$ is odd. Therefore,
$$
        N_{\cO_k}(n,K,m)= \begin{cases}
            S_{\cO_K}(n)\cdot q^{nm}+E_{n,m} &\textrm{if $m\ge 2g_K$},\\
       0 &\textrm{if $m=2g_K-1$}.
        \end{cases}
        $$
\vskip 1mm
\noindent
$\bullet$ Assume finally that $\infty$ remains inert in $K/k$, and let $P$ be the unique place of $K/k$ lying above $\infty$. In particular, $\deg(P)=2$ and $U_K=2\N$. Note that $$
N_{\cO_k}(n,K,m)=\#(\sL_n'(\frac{m}{2}P) \setminus \cO_k^n).
$$
If $m-1 \geq g_K \geq 1 + \frac{m}{2}$, Lemma \ref{lem:inert_rami} again implies that
$$
N_{\cO_k}(n,K,m)=0.
$$
Suppose now that $g_K < 1 + \frac{m}{2}$. By Lemma \ref{clm:inert1}, 
$$
 \#\left\lbrace\alpha\in \cO_K^n, (\alpha)_\infty=\frac{m}{2}P\right\rbrace=
 \begin{cases}
      S_{\cO_K}(n)\cdot q^{nm} &\textrm{if $m\ge 2g_K+2$},\\
        \frac{q^n}{q^n+1}S_{\cO_K}(n)\cdot q^{nm} &\textrm{ if $m=2g_K$}.
 \end{cases}
$$
As in the ramified case, there are $O(q^\frac{mn}{2})$ elements  $\alpha\in \cO_k^n$ with $(\alpha)_\infty=\frac{m}{2}P$. Thus, 
$$
        N_{\cO_k}(n,K,m)= \begin{cases}
            S_{\cO_K}(n)\cdot q^{nm}+E_{n,m} &\textrm{if $m\ge 2g_K+2$},\\
        \frac{q^n}{q^n+1}S_{\cO_K}(n)\cdot q^{nm}+E_{n,m} &\textrm{if $m=2g_K$}.
        \end{cases}
        $$ 
\end{proof} 

\begin{lem}\label{high_gen_count}
Assume $q$ is odd. For $m>4$, we have
    \begin{equation}\label{eq:sumsympttt}
        \sum\limits_{\substack{[K:k]=2\\q_K=q\\\frac
        {m}{2}+1\le g_K}} N_{\cO_k}(3,K,m)\ll mq^{\frac{5m}{2}}.
    \end{equation}
\end{lem}
\begin{proof}
    By Lemma \ref{lem:n=2_fixed_K_count}, 

   \begin{align*}
       \sum\limits_{\substack{[K:k]=2\\q_K=q\\\frac
        m2+1\le g_K}} N_{\cO_k}(3,K,m)&= \sum\limits_{\substack{[K:k]=2\\q_K=q\\\frac
        m2+1\le g_K\le m-1}} N_{\cO_k}(3,K,m)\\&=
        \sum\limits_{i=\lceil\frac{m}{2}\rceil+1}^{m-1}\;\;\sum\limits_{\substack{[K:k]=2\\q_K=q,s_K=2\\g_K=i}} 2\left(m-g_K+\frac{1}{q^3-1}\right)\left(1-\frac{1}{q^3}\right)^2q^{3(m+1-g_K)}\\
        &\le \sum\limits_{i=\lceil\frac{m}{2}\rceil+1}^{m-1}\;\;\sum\limits_{\substack{[K:k]=2\\q_K=q,s_K=2\\g_K=i}} m q^{3(m+1-g_K)}.\\
   \end{align*}
Consequently, Lemma \ref{lem:gen_count} gives
\begin{alignat*}3
       \sum\limits_{\substack{[K:k]=2\\q_k=q\\\frac
        m2+1\le g_K}} N_{\cO_k}(3,K,m) \ll \sum\limits_{i=\lceil\frac{m}{2}\rceil+1}^{m-1} m q^{3(m+1)} q^{-i}\ll mq^{\frac{5m}{2}}.
   \end{alignat*}
This completes the proof.
\end{proof}

\subsection{Summation over quadratic extensions}

\begin{proof}[Proof of Theorem \ref{thm:main2} in the case $(n,e)=(3,2)$]
By Corollary 1.4 and Lemma \ref{high_gen_count}, 

\begin{align*}
    N_{\cO_k}(3,2,m)&= \sum \limits_{\substack{[K:k]=2\\m> 2g_K-2}}N_{\cO_k}(3,K,m)+ \sum \limits_{\substack{[K:k]=2\\m\le 2g_K-2}}N_{\cO_k}(3,K,m)\\
    &= \sum \limits_{\substack{[K:k]=2\\m> 2g_K-2}} \left[S_{\cO_K}(3) m^{s_K-1}\, q^{3m}
+ O_{k}(\tilde E(m))\right]+ O_k(mq^{\frac{5m}{2}}),
\end{align*}
where
$$
\tilde E(m)=
\begin{cases}
m^{s_K-2} q^{3m-3g_K} + q^{2m} &\textrm{if } s_K \geq 2,\\
q^{2m}&\textrm{if } s_K=1.
\end{cases}
$$
A direct computation using Lemma \ref{lem:gen_count}, together with the fact that $s_K\le 2$, shows that

$$\sum \limits_{\substack{[K:k]=2\\m> 2g_K-2}}  \tilde E(m)=O_k(q^{3m}).$$
This completes the proof of the theorem in the case $(n,e)=(3,2)$.
\end{proof}

\begin{proof}[Proof of Theorem \ref{thm:main2} in the case $(n,e)=(2,2)$] To obtain the asymptotic estimate for $N_{\cO_k}(2,2,m)$, we use Lemma \ref{lem:n=2_fixed_K_count}. We have
\begin{alignat*}1
    N_{\cO_k}(2,2,m)&=\sum\limits_{\substack{[K:k]=2\\q_K=q}} N_{\cO_k}(2,K,m)\\
    &=\sum\limits_{\substack{[K:k]=2\\s_K=2\\q_K=q}} N_{\cO_k}(2,K,m)+ \sum\limits_{\substack{[K:k]=2\\s_K=1\\q_K=q}} N_{\cO_k}(2,K,m).\\
\end{alignat*}

Using Lemma \ref{lem:n=2_fixed_K_count}, we see that 
\begin{align*}
\sum\limits_{\substack{[K:k]=2\\s_K=1\\q_K=q}} N_{\cO_k}(2,K,m)&\ll \sum\limits_{\substack{[K:k]=2\\g_K<1+m/2\\q_K=q}}q^{2(m-g_K)}+ \sum\limits_{\substack{[K:k]=2\\g_K< 1+m/2\\q_K=q}}E_m\\
    &\ll \sum\limits_{i=0}^{1+\frac{m}{2}}q^{2(m-i)}q^{2i}+\sum\limits_{i=0}^{1+\frac{m}{2}}E_m q^{2i} &\textrm{ by Lemma \ref{lem:gen_count}}\\
&\ll mq^{2m}.
\end{align*}

We may therefore focus on the first term. Using Lemma \ref{lem:n=2_fixed_K_count}, we see that
\begin{alignat*}1
    \sum\limits_{\substack{[K:k]=2\\s_K=2\\q_K=q}} N_{\cO_k}(2,K,m)=\sum\limits_{\substack{[K:k]=2\\s_K=2\\q_K=q\\1+\frac{m}{2}\le g_K\le m-1}}2(m-g_K)S_{\cO_K}(2)\cdot q^{2m}+ \sum\limits_{\substack{[K:k]=2\\s_K=2\\q_K=q\\g_K<1+\frac{m}{2}}} m S_{\cO_K}(2)\cdot q^{2m}+\\
    \\
    O\left(\sum\limits_{\substack{[K:k]=2\\s_K=2\\q_K=q\\1+\frac{m}{2}\le g_K\le m-1}}S_{\cO_K}(2)q^{2m}+ \left[\sum\limits_{\substack{[K:k]=2\\s_K=2\\q_K=q\\g_K<1+\frac{m}{2}}}S_{\cO_K}(2) q^{2m}+E_m\right]\right).
\end{alignat*}

 Using Lemma \ref{lem:gen_count}, we obtain that 
 \begin{align*}
     \sum\limits_{\substack{[K:k]=2\\s_K=2\\q_K=q\\1+\frac{m}{2}\le g_K\le m-1}}S_{\cO_K}(2)q^{2m}+ \left[\sum\limits_{\substack{[K:k]=2\\s_K=
     2\\q_K=q\\g_K<1+\frac{m}{2}}}S_{\cO_K}(2) q^{2m}+E_m\right]&\ll 
     \sum\limits_{i=\lceil1+\frac{m}{2}\rceil}^{m-1}q^{2m}+ \left[\sum\limits_{i=0}^{\lfloor 1+\frac{m}{2}\rfloor} q^{2m}+q^{m+2i}\right]\\
     &
     \ll mq^{2m}.
 \end{align*}
Thus
 $$N_{\cO_k}(2,2,m)=q^{2m}\left[\sum\limits_{\substack{[K:k]=2\\s_K=2\\q_K=q\\g_K<1+\frac{m}{2}}}mS_{\cO_K}(2)+\sum\limits_{\substack{[K:k]=2\\s_K=2\\q_K=q\\1+\frac{m}{2}\le g_K<m}}2(m-g_K)S_{\cO_K}(2)\right]+O(mq^{2m}).$$

Using the precise count in Lemma \ref{lem:gen_count}, we now obtain

$$
N_{\cO_k}(2,2,m)=\frac{3}{4}q^3(q-1)\left(1-\frac{1}{q^2} \right)^2 m^2 q^{2m}+O(m q^{2m}).
$$
This completes the proof.
\end{proof}
\begin{proof}[Proof of Theorem \ref{thm:main2} in the case $(n,e)=(1,2)$]
We note that every \emph{irreducible} polynomial of the form $f=X^2+AX+B\in \cO_k[X]$, with $h(f)=m$ uniquely corresponds to a pair of elements 
$(\alpha,\bar\alpha)$ (which are conjuagtes), given by the roots of $f$. It is clear that $k(\alpha)=k(\bar\alpha)$, $[k(\alpha):k]=2$ and  $h(\alpha)=h(\bar\alpha)$, since $\alpha,\bar\alpha$ are conjugates. By Lemma \ref{lem:heighaddpoly}, it follows that $h(\alpha)=h(\bar\alpha)=\frac{m}{2}.$

We therefore have using Lemma \ref{lem:heighaddpoly}
\begin{multline*}
N_{\cO_k}(1,2,m)+\sum\limits_{i=0}^{m}N_{\cO_k}(1,1,i)\cdot N_{\cO_k}(1,1,m-i)+N_{\cO_k}(1,1,\frac{m}{2})\\
= 2\times \#\{f=X^2+AX+B\in \cO_k[X], h(f)=m\}.
\end{multline*}
It is easy to see that 
$$
\#\{f=X^2+AX+B\in \cO_k[X], h(f)=m\}= (q^2-1)q^{2m},
$$
and
$$
N_{\cO_k}(1,1,i)=\begin{cases}
(q-1)q^i &\textrm{ if } i>0\\
q &\textrm{ if } i=0.
\end{cases}
$$
Thus,
$$
\sum\limits_{i=0}^{m}N_{\cO_k}(1,1,i)\cdot N_{\cO_k}(1,1,m-i)=(m-1)(q-1)^2q^m+2(q-1)q^{m+1}.
$$
Therefore,
$$
N_{\cO_k}(1,2,m)=2\left(q^2-1\right)q^{2m}+O(mq^{m}).
$$
This completes the proof
\end{proof}

\section{Proof of Theorem \ref{hom_forms_gen} and Corollary \ref{cor: hom_forms_gen}}\label{sec:last}

Throughout this section, $d=[k:\F_q(t)]$ and $p={\rm char}(k)$. For $e \geq 1$, define
$$
N^{\textrm{sep}}_{\cO_k}(n,e,m)=\sum\limits_{\substack{[K:k]=e\\q_K=q\\ K/k\,{\rm separable}}}N_{\cO_k}(n,K,m).
$$ Similarly, let $N_{\cO_k}^\textrm{sep}(n,d,m,\textrm{hom})$ count the forms in $N_{\cO_k}(n,e,m,\textrm{hom})$ whose linear factors (as in \eqref{eq:decomlin}) generate separable extensions of $k$. For any form as in \eqref{eq:decomlin}, we call \emph{the factor field} $F$ the field generated by one of its factor over $k$.

\subsection{Some lemmas} Let $K/k$ be a finite extension of degree $e$ with $q_K=q$, and let $K_s$ be the maximal separable extension of $K/k$. Write
$$
[K:K_s]=p^\ell.
$$
Since $k$ is a one-variable function field over $\F_q$, one has
$$
K_s=K^{p^\ell}.
$$

\begin{lem}\label{lem:overKK}
Let $K/k$ and $\ell$ be as above, and suppose that $\ell \geq 1$. Then
\begin{equation}\label{eq:noknnn}
N_{\cO_k}(n,K,m)=N_{\cO_k}(n,K_s,m)-N_{\cO_k}(n,K_s,m/p).
\end{equation}
\end{lem}

\begin{proof}
The map 
$$\Phi_\ell:\alpha \in K^n \mapsto \alpha^{p^\ell} \in K_s^n
$$ is a bijection. Let $P$ be a place of $K_s$, and let $Q$ be its unique extension to $K$. Since $K/K_s$ is purely inseparable of degree $p^\ell$, the place $Q$ is totally ramified. For $x \in K$, we have
$$
\ord_Q(x)=\ord_P(x^{p^\ell}).
$$
Consequently, $\alpha \in \cO_{K}^n$ if and only if $\Phi_\ell(\alpha) \in \cO_{K_s}^n$, and
    $$h_K(\alpha)=h_K([1:\alpha_1:\dots: \alpha_n])= h_{K_s}([1:\alpha_1^{p^\ell}:\dots: \alpha_n^{p^\ell}])=h_{K_s}(\Phi_\ell(\alpha)).
    $$ 

Now let $\beta=\Phi_\ell(\alpha)$. If $k(\alpha)=K$, then $k(\beta)=K_s$. Assume $k(\beta)=K_s$. We first note that, for $x \in K$,
$$
x\in K^p \iff x^{p^{\ell-1}} \in K_s.
$$
Thus, the point $\alpha$ does not generate $K$ over $k$ precisely when $\beta^{1/p} \in K_s^n$. Therefore, the excluded vectors from the ones counted $N_{\cO_k}(n,K_s,m)$ are exactly exactly $\beta=\gamma^p$ such that $k(\gamma^p)=K_s$--hence $k(\gamma)=K_s$ as $K_s/k$ is separable--and $h_{K_s}(\gamma)=h_{K_s}(\beta)/p$. 
        
   Consequently,
   $$
   N_{\cO_k}(n,K,m)=N_{\cO_k}(n,K_s,m)-N_{\cO_k}(n,K_s,m/p).
   $$
\end{proof}

For $\ell=0,\dots,v_p(e)$, let $N^{(\ell)}_{\cO_k}(n,e,m)$ count the points counted by $N_{\cO_k}(n,e,m)$ whose generated extension has inseparable degree $p^\ell$ over $k$.

\begin{lem}
Let $r=v_p(e)$. Then
$$
N^{(0)}_{\cO_k}(n,e,m)=N^{{\rm sep}}_{\cO_k}(n,e,m),
$$
and, for $\ell=1,\dots,r$,
\begin{equation}\label{eq:nnmin}
N^{(\ell)}_{\cO_k}(n,e,m)=N^{{\rm sep}}_{\cO_k}(n,e/p^\ell,m)-N^{{\rm sep}}_{\cO_k}(n,e/p^\ell,m/p).
\end{equation}
Consequently,
\begin{equation}\label{eq:susmsusm}
N_{\cO_k}(n,e,m)=N_{\cO_k}^{\rm{sep}}(n,e,m)+\sum_{\ell=1}^r\left(N_{\cO_k}^{{\rm sep}}(n,e/p^\ell,m)-N_{\cO_k}^{{\rm sep}}(n,e/p^\ell,m/p)\right).
\end{equation}
\end{lem}

\begin{proof}
For a fixed $\ell$, the map $K \mapsto K_s$ is a bijection preserving the constant field between degree-$e$ extensions of inseparable degree $p^\ell$ and separable extensions of degree $e/p^\ell$. Summing \eqref{eq:noknnn} over the corresponding fields gives \eqref{eq:nnmin}, and summing over $\ell$ gives \eqref{eq:susmsusm}.
\end{proof}

\begin{cor}\label{exp_sep}\label{cor: exp_sep} Let $r=v_p(e)$.
We have 
$$
N_{\cO_k}^{{\rm sep}}(n,e,m)=
\begin{cases}
 N_{\cO_k}(n,e,m)-N_{\cO_k}(n,e/p,m)+N_{\cO_k}^{\rm{sep}}(n,e/p,m/p) &\textrm{ if } p\mid e, \\
 N_{\cO_k}(n,e,m) &\textrm{ otherwise.}
 \end{cases}
    $$
Iterating the formula above yields 
    \begin{equation}\label{eq:iterzeqa}
    N_{\cO_k}^{\rm{sep}}(n,e,m)=\sum\limits_{i=0}^{r}\left[ N_{\cO_k}(n,e/p^i,m/p^i)-N_{\cO_k}(n,e/p^{i+1},m/p^i)\right].
    \end{equation}
\end{cor}
\begin{proof}
    The claim is immediate when $p\nmid e$. Suppose therefore that $p \mid e$. Subtracting the decomposition \eqref{eq:susmsusm} for degree $e/p$ from the decomposition for degree $e$ gives
    $$
N_{\cO_k}(n,e,m)-N_{\cO_k}(n,e/p,m)=N_{\cO_k}^{{\rm sep}}(n,e,m)-N_{\cO_k}^{\rm{sep}}(n,e/p,m/p).
$$
\end{proof}

We now establish the analogue of the above decomposition for forms. Suppose that a form $F$ counted by $N_{\cO_k}(n,e,m,{\rm hom})$ has a factor field $K$ of inseparable degree $p^\ell$. Thus it has $e/p^\ell$ distinct linear factors, each occurring with multiplicity $p^\ell$. After re-indexing the distinct factor as
$$
L_i=x_0+\sum_{j=1}^n\alpha_{i,j}x_j \quad \left(i=1,\dots, \frac{e}{p^\ell} \right),
$$ set
$$
L'_i=x_0+\sum_{i=1}^n \alpha_{i,j}^{p^\ell}x_j \in K_s[x_0,\dots,x_n] \textrm{ and }F'=\prod_{i=1}^{e/p^\ell} L'_i \in k[x_0,\dots,x_n].
$$
Then
\begin{equation}\label{eq:foffor}
F(x_0,\dots,x_n)=F'(x_0^{p^\ell},\dots,x_n^{p^\ell}).
\end{equation}
In particular, $F$ and $F'$ have the same nonzero coefficients, so
$$
h(F)=h(F').
$$
The form $F'$ has degree $e/p^\ell$, and its factor field is the separable field $K_s$. Moreover, $F'$ is irreducible over $k$. 
\begin{lem}\label{sep_count}\label{cor: exp_sephom}
    We have 
    $$ N_{\cO_k}(n,e,m,{\rm hom})=N_{\cO_k}^{\rm{sep}}(n,e,m,{\rm hom})+\sum_{\ell=1}^r\left(N_{\cO_k}^{{\rm sep}}(n,e/p^\ell,m,{\rm hom})-N_{\cO_k}^{{\rm sep}}(n,e/p^\ell,m/p,{\rm hom})\right).
    $$
\end{lem}

\begin{proof}
 By the argument above, the assignment $F \mapsto F'$ is an injective map between decomposable homogeneous forms of degree $e$ whose factor field  has maximal separable extension of degree $e/p^\ell$ over $k$, and decomposable homogeneous forms of degree $e/p^\ell$ whose factor field is separable extension of degree $e/p^\ell$ over $k$.  Using the same argument as for $N_{\cO_k}(n,e,m)$ completes the proof.
\end{proof}

\begin{lem}\label{hom_aff_comp} For every $e\geq 1$,
    $$
    eN_{\cO_k}^{{\rm sep}}(n,e,m,{{\rm hom}})= N_{\cO_k}^{{\rm sep}}(n,e,m).
    $$
\end{lem}

\begin{proof}
  Let $F$ be counted by $N_{\cO_k}^{{\rm sep}}(n,e,m,{{\rm hom}})$. Since its factor field is separable, the factorization \eqref{eq:decomlin} has $e$ distinct conjugate factors.
 By Lemma \ref{lem:heighaddpoly},   
  $$
  h(F)= \sum_{j=1}^e h(L_j)= e
  h(L_i), \quad i=1, \dots,e.
  $$ Since $h(F)=m/d$, the coefficient vector of each factor has height $m/(ed)$. Moreover, each coefficient is integral over $\cO_k$: for every $j$, the specialization
  $$
F(x_0,0,\dots,0,x_j,0,\dots,0)=\prod_{i=1}^e(x_0+\alpha_{i,j}x_j)
  $$ is monic in $x_0$ and has coefficients in $\cO_k$. Setting $x_j=1$ shows that $-\alpha_{i,j}$ is a root of a monic polynomial in $\cO_k[x_0]$.

  Conversely, let $\alpha \in \cO_{\overline{k}}^n$ be counted by $N_{\cO_k}^{{\rm sep}}(n,e,m)$, and set
  $$
L=x_0+\alpha_{1}x_1+\dots+\alpha_nx_n.
  $$
  The norm form
  $$
 F=\prod_{\sigma}\sigma(L)
  $$
  where $\sigma$ ranges over the $k$-embeddings of $k(\alpha)$ into $\overline{k}$, belongs to $\cO_k[x_0,\dots,x_n]$. Moreover, $F$ is irreducible over $k$, has height $m/d$ and such that its factor field has constant field $\F_q$.

Therefore each form $F$ is exactly a norm from of $e$ such $L$, namely the coefficient vectors of its $e$ conjugate factors.
  
Consequently,
$$
eN_{\cO_k}^{{\rm sep}}(n,e,m,{{\rm hom}})= N_{\cO_k}^{{\rm sep}}(n,e,m),
$$ as claimed.
  \end{proof}

\subsection{Proof of Theorem \ref{hom_forms_gen}} By Lemma \ref{cor: exp_sephom} and Lemma \ref{hom_aff_comp}, 

 \begin{align*}
     N_{\cO_k}(n,e,m,\text{hom})&=N_{\cO_k}^{\rm sep}(n,e,m,{\rm hom})+\sum\limits_{\ell=1}^{r}\left[N_{\cO_k}^{\textrm{sep}}(n,e/p^\ell,m,\text{hom})-N_{\cO_k}^{\textrm{sep}}(n,e/p^\ell,m/p,{\rm hom})\right]\\
     &=\frac{1}{e}N_{\cO_k}^{\rm sep}(n,e,m)+\sum\limits_{\ell=1}^{r}\frac{p^\ell}{e}\left[N_{\cO_k}^{{\rm sep}}(n,e/p^\ell,m)-N_{\cO_k}^{\textrm{sep}}(n,e/p^\ell,m/p)\right].
\end{align*}
From Corollary \ref{exp_sep}, it follows that
\begin{align*}
   N_{\cO_k}(n,e,m,\text{hom}) & =\frac{1}{e}\left[N_{\cO_k}(n,e,m)-N_{\cO_k}(n,e/p,m)+N_{\cO_k}^{\rm sep}(n,e/p,m/p)\right]\\
    &+\sum\limits_{\ell=1}^{r}\frac{p^\ell}{e}\left(N_{\cO_k}(n,e/p^\ell,m)-N_{\cO_k}(n,e/p^{\ell+1},m)\right)\\
    &-\sum\limits_{\ell=1}^{r}\frac{p^\ell}{e}\left(N_{\cO_k}^{\rm sep}(n,e/p^\ell,m/p)-N_{\cO_k}^{\rm sep}(n,e/p^{\ell+1},m/p)\right)\\
    &= \frac{1}{e}N_{\cO_k}(n,e,m)+\frac{1}{e}\sum\limits_{\ell=1}^{r}(p^\ell-p^{\ell-1})\left(N_{\cO_k}(n,e/p^\ell,m)-N_{\cO_k}^{\rm sep}(n,e/p^\ell,m/p)\right).
\end{align*}
By \eqref{eq:iterzeqa},
\begin{align*}
   N_{\cO_k}(n,e,m,\text{hom}) 
     &= \frac{1}{e}N_{\cO_k}(n,e,m)+\frac{1}{e}\sum\limits_{\ell=1}^{r}(p^\ell-p^{\ell-1}) N_{\cO_k}(n,e/p^\ell,m)\\
    &-\frac{1}{e}\sum\limits_{\ell=1}^{r}(p^\ell-p^{\ell-1})\sum\limits_{j=\ell}^{r}\left(N_{\cO_k}(n,e/p^j,m/p^{j-\ell+1})-N_{\cO_k}(n,e/p^{j+1},m/p^{j-\ell+1})\right).
 \end{align*}
Letting $j'=j-\ell$, we have
\begin{align*}
  N_{\cO_k}(n,e,m,\rm{ hom})=
     &\frac{1}{e}N_{\cO_k}(n,e,m)+\frac{p-1}{e}\sum\limits_{\ell=1}^{r}p^{\ell-1} N_{\cO_k}(n,e/p^\ell,m)\\
     &-\frac{p-1}{e}\sum\limits_{\ell=1}^{r}p^{\ell-1}\sum\limits_{j'=0}^{r-\ell}\left(N_{\cO_k}(n,e/p^
     {j'+\ell},m/p^{j'+1})-N_{\cO_k}(n,e/p^{j'+\ell+1},m/p^{j'+1})\right).\\
    \end{align*}
On the other hand,

    \begin{align*}
        \sum\limits_{\ell=1}^{r}p^{\ell-1}\sum\limits_{j'=0}^{r-\ell}\left(N_{\cO_k}(n,e/p^
     {j'+\ell},m/p^{j'+1})-N_{\cO_k}(n,e/p^{j'+\ell+1},m/p^{j'+1})\right)\\ 
      =\sum\limits_{j'=0}^{r-1}\sum\limits_{\ell=1}^{r-j'}p^{\ell-1}\left(N_{\cO_k}(n,e/p^
     {j'+\ell},m/p^{j'+1})-N_{\cO_k}(n,e/p^{j'+\ell+1},m/p^{j'+1})\right)\\
     =\sum\limits_{j'=0}^{r-2}\sum\limits_{\ell=2}^{r-j'}(p^{\ell-1}-p^{\ell-2})N_{\cO_k}(n,e/p^{\ell+j'},m/p^{j'+1})+\sum\limits_{j'=0}^{r-1}N_{\cO_k}(n,e/p^{j'+1},m/p^{j'+1})\\
     =(p-1)\sum\limits_{j'=0}^{r-2}\sum\limits_{\ell=2}^{r-j'}p^{\ell-2}N_{\cO_k}(n,e/p^{\ell+j'},m/p^{j'+1})+\sum\limits_{j'=0}^{r-1}N_{\cO_k}(n,e/p^{j'+1},m/p^{j'+1}).
    \end{align*}
 Therefore,

    \begin{align*}
          N_{\cO_k}(n,e,m,\rm{ hom})=
     &\frac{1}{e}N_{\cO_k}(n,e,m)+\frac{p-1}{e}\sum\limits_{\ell=1}^{r}p^{\ell-1} N_{\cO_k}(n,e/p^\ell,m)\\
     &-\frac{p-1}{e}\sum\limits_{j'=0}^{r-1}N_{\cO_k}(n,e/p^{j'+1},m/p^{j'+1})\\
     &-\frac{(p-1)^2}{e}\sum\limits_{j'=0}^{r-2}\sum\limits_{\ell=2}^{r-j'}p^{\ell-2}N_{\cO_k}(n,e/p^{\ell+j'},m/p^{j'+1}).
    \end{align*}
Recalling that
   $$ T_i=N_{\cO_k}\left(n,\frac{e}{p^i},\frac{m}{p^i}\right)+(p-1)\sum\limits_{j=i+1}^{r}p^{j-i-1}N_{\cO_k}\left(n,\frac{e}{p^j},\frac{m}{p^i}\right),$$
we obtain 
   $$N_{\cO_k}(n,e,m,\rm{ hom})=\frac{1}{e}\left(T_0-(p-1)\sum\limits_{i=1}^{r}T_i\right).$$
 This completes the proof of Theorem \ref{hom_forms_gen}.

The proof of Corollary \ref{cor: hom_forms_gen} is straightforward- it follows directly from Theorem \ref{hom_forms_gen}, Theorem 
 \ref{thm:main1} and Theorem \ref{thm:main2}.
\vskip 1mm
\noindent
{\bf Acknowledgements}: The authors would like to express their sincere gratitude to Martin Widmer for suggesting this project and for his careful review of the manuscript. His insightful comments and constructive suggestions significantly contributed to improving the quality of this work. The first author is supported by a JSPS Postdoctoral Fellowship, and also acknowledges Noy Soffer Aranov, Christian Elsholtz and TU Graz for funding his research visits at TU Graz between March and June 2026.
\vskip 1mm
\noindent
{\bf Use of AI}: ChatGPT was used solely for grammatical corrections and for checking the lemmas in Subsection 6.1. All the ideas and arguments in this paper are due to the authors, who take full responsibility for them.

\bibliography{Biblio2}
\bibliographystyle{alpha}
\end{document}